\documentclass[11pt]{amsart}

\usepackage{amsmath, amssymb, amsfonts, amsthm}
\usepackage{graphicx}
\usepackage{color}
\usepackage{bm}
\usepackage{subfigure}
\usepackage{graphicx}
\usepackage{mathabx}
\usepackage{multirow}
\usepackage{setspace}
\usepackage[title]{appendix}
\usepackage{mathrsfs} 
\usepackage{enumerate}
\usepackage{comment}
\usepackage{longtable}
\usepackage{tikz}
\usepackage{wrapfig}
\usepackage{float}

\usetikzlibrary{matrix,arrows,decorations.pathmorphing}						
							
\usepackage{amsthm}
\theoremstyle{plain}
\newtheorem{theorem}{Theorem}[section]

\theoremstyle{plain}
\newtheorem{lemma}[theorem]{Lemma}

\theoremstyle{plain}
\newtheorem{Prop}{Proposition}[section]

\theoremstyle{plain}
\newtheorem{corollary}[theorem]{Corollary}

\theoremstyle{plain}

\theoremstyle{definition}
\newtheorem{Def}{Definition}[section]
\theoremstyle{remark} 
\newtheorem{remark}{Remark}

\theoremstyle{definition}

\theoremstyle{definition} 
\newtheorem{Example}{Example}

\def\ed{\mathrm{d}}

\def\I{\mathcal I}

\def\J{\mathcal J}

\def\R{\mathbb R}

\def\P{\mathbb P}
\def\S{\mathbb S}

\def\E{\mathbb E}

\def\e{{\bf e}}
\def\x{{\bf x}}
\def\u{{\bf u}}

\def\diag{{\rm diag}}

\def\W{\wedge}
\def\<{\langle}
\def\>{\rangle}

\def\SO{{\rm SO}}

\def\GL{{\rm GL}}

\def\SL{{\rm SL}}

\def\alg{{\rm alg}}

\def\tT{\tilde T}

\def\({\left(}
\def\){\right)}

\def\lbb{[\![}
\def\rbb{]\!]}

\def\bs{\boldsymbol}

\def\mw{\&$\W$}

\def\codesp{3.4mm}

\DeclareFontFamily{U}{MnSymbolC}{}
\DeclareSymbolFont{MnSyC}{U}{MnSymbolC}{m}{n}
\DeclareFontShape{U}{MnSymbolC}{m}{n}{
    <-6>  MnSymbolC5
   <6-7>  MnSymbolC6
   <7-8>  MnSymbolC7
   <8-9>  MnSymbolC8
   <9-10> MnSymbolC9
  <10-12> MnSymbolC10
  <12->   MnSymbolC12}{}
\DeclareMathSymbol{\lefthook}{\mathbin}{MnSyC}{'270}

\title[Geometry of B\"acklund Transformations II]{Geometry of B\"acklund Transformations II: Monge-Amp\`ere Invariants}
\author{Yuhao Hu} 
\address{Department of Mathematics, 395 UCB, University of
Colorado, Boulder, CO 80309-0395}
\email{Yuhao.Hu@colorado.edu}


\begin{document}

\maketitle
\begin{abstract}
	This paper is concerned with the question:
	\emph{For which pairs of hyperbolic Euler-Lagrange systems in the plane does there exist
	a rank-$1$ B\"acklund transformation relating them?} We express some obstructions
	to such existence in terms
	of the local invariants of the Euler-Lagrange systems.
	In addition, we discover a class of B\"acklund transformations relating two hyperbolic Euler-Lagrange
	systems of distinct types.
\end{abstract}

\setcounter{tocdepth}{1}
\tableofcontents

\section{Introduction}
			
A B\"acklund transformation is a way to relate solutions of two PDE systems $\mathcal{E}_1$ and
$\mathcal{E}_2$, in such a manner that, given a solution of $\mathcal{E}_1$, one can use it to obtain
solutions of $\mathcal{E}_2$ by solving ODEs, and \emph{vice versa}. Early studies of 
B\"acklund transformations date back to the late 19th century.

\subsection{Three Examples}
The following examples of B\"acklund transformations are classical.
\vskip  2mm
{\bf 1.} Let $u = u(x,y)$ be a harmonic function. One can find a $1$-parameter family of harmonic functions
		$v = v(x,y)$ by substituing $u(x,y)$ into the Cauchy-Riemann system in $u$ and $v$ and by
		solving ODEs. In this sense, the Cauchy-Riemann system is a B\"acklund transformation relating
		solutions of the Laplace equation $\Delta u(x,y) = 0$.
\vskip 1mm
{\bf 2.} Let $u = u(x,y)$ be a solution of the \emph{sine-Gordon} equation 
		\begin{equation}
			u_{xy} = \frac{1}{2} \sin(2u).	\label{SG}
		\end{equation}
	Substituting $u(x,y)$ in the system (with parameter $\lambda \ne 0$) in $u$ and $v$:
	\begin{equation}
		\left\{
			\begin{array}{l}
				u_x - v_x = \lambda \sin(u+v),\\[0.8em]
				u_y+v_y = \lambda^{-1}\sin(u-v),
			\end{array}
		\right.		\label{SGT}
	\end{equation}
	one obtains a compatible first-order PDE system in $v(x,y)$, which can be solved using ODE methods.
	In this case, the system \eqref{SGT} is a B\"acklund transformation relating solutions of 
	the sine-Gordon equation. In particular, if we start by setting $u(x,y) = 0$, which is a solution of 
	\eqref{SG}, solving \eqref{SGT} (with a fixed $\lambda$) for $v(x,y)$ yields a 1-parameter family 
	of solutions
	of \eqref{SG}:
			\[
				v(x,y) = \arctan\(C \exp\(-\lambda x - \lambda^{-1}y\)\),
			\]
	where $C$ is an arbitrary constant.		
	Such $v(x,y)$ are known as \emph{$1$-soliton} solutions of \eqref{SG}.
	Iterative application of the B\"acklund transformation \eqref{SGT} yields the so-called
	\emph{$n$-solitons} of \eqref{SG}. 
	Interested reader can see \cite{terng2000geometry} and \cite{Rogers} for 
	more details.
\vskip 1mm
{\bf 3.} Given two immersed surfaces $S_1, S_2\subset \E^3$, by a \emph{pseudo-spherical (p.s.) line
congruence between them} we mean immersions $\sigma_1,\sigma_2: U\rightarrow \E^3$ ($U\subset \R^2$ open) satisfying
$\sigma_i(U) = S_i$ and
	\begin{quote}
	 \begin{enumerate}[{$(1)$}]
	\item{The distance $d_{\E^3}(\sigma_1(p), \sigma_2(p))$ is a constant $r>0$;}
	\item{For each $p\in U$, the line through $\sigma_1(p)$ and $\sigma_2(p)$ is tangent to both surfaces
			at $\sigma_1(p)$ and $\sigma_2(p)$, respectively;}
		\item{The respective normals $\bs n_1(p)$ and $\bs n_2(p)$ form a constant angle $\theta\in (0,\pi)$.}
	 \end{enumerate}
	 \end{quote}
	When these conditions hold, $S_1, S_2$ are called the \emph{focal surfaces} of the corresponding 
	p.s. line congruence.
	
	 It is a theorem of Bianchi that
	\emph{$S_1, S_2\subset \E^3$ admit a p.s. line conguence (with parameters $r,\theta$) 
	between them only if
	 they both have the Gauss curvature $K = -\sin^2\theta/r^2$. Conversely, 
	 if $S\subset \E^3$ satisfies $K = -\sin^2\theta/r^2$ for some constants $r>0$ and 
	 $\theta\in (0,\pi)$,
	 then one can construct
	 a 1-parameter family of p.s. line congruences with parameters $r$ and $\theta$ and
	 with $S$ being a focal surface.} 
	 
	 In fact, when $S_1,S_2$ admit a p.s. line congruence with parameters $r$ and $\theta$, 
	 the Gauss and Codazzi equations of $S_i$ $(i = 1,2)$ together 
	 form a PDE system $\mathcal{E}_{r, \theta}$. The integrability 
	 of $\mathcal{E}_{r,\theta}$ 
	 puts conditions on the respective Gauss
	 curvatures.
	 On the other hand, given
	 a surface $S$ with constant Gauss curvature $K =  -\sin^2\theta/r^2 <0$, 
	 the problem of finding a p.s.
	 line congruence with parameters $r,\theta$ and with $S$ as a focal surface
	 reduces to integrating a Frobenius (aka. completely integrable) system.
	 In this sense, $\mathcal{E}_{r,\theta}$ is a B\"acklund transformation relating surfaces
	 with $K =  -\sin^2\theta/r^2$. For more details, including
	 how $\mathcal{E}_{r,\theta}$ relates to the sine-Gordon equation
	 and its B\"acklund transformations,  see \cite{chern-terng} and \cite{BGG}.\\

	Using the theory of exterior differential systems, one can study B\"acklund transformations
	from a geometric viewpoint. For instance,
	the B\"acklund transformation $\mathcal{E}_{r,\theta}$ in Example {\bf 3} 
	 above is \emph{homogeneous} in the sense that
	  the symmetry group of the system acts locally transitively on the space of variables.
	 It has \emph{rank} $1$ in the sense that
	 it generates
	 a 1-parameter family of surfaces with $K =  -\sin^2\theta/r^2$ from a given one.
	 The PDE system for $S$ with a prescribed constant Gauss curvature $K<0$
	  is an example of a \emph{hyperbolic Monge-Amp\`ere system}.
	 A complete classification of rank-$1$ homogeneous B\"acklund transformations relating two hyperbolic 
	 Monge-Amp\`ere systems has been obtained by Jeanne N. Clelland in \cite{C01}, using Cartan's method
	 of equivalence.
	 
\subsection{Geometric Formulations}	 

Classically, a \emph{Monge-Amp\`ere equation} in $z(x,y)$ is a second-order PDE of the form
	\begin{equation}
		A(z_{xx}z_{yy} - z_{xy}^2)+Bz_{xx}+2Cz_{xy}+Dz_{yy}+E = 0,	\label{MAequation}
	\end{equation}
where $A,B,C,D,E$ are functions of $x,y,z,z_x,z_y$. This is a class of 
 PDEs that are closely related to surface geometry and calculus of variations (see 
\cite{bryantedslectures} and \cite{BGG}).
A Monge-Amp\`ere equation \eqref{MAequation} is said to be
\emph{hyperbolic} (resp., \emph{elliptic}, \emph{parabolic}) if $AE - BD+ C^2$ is positive (resp.,  negative, zero).

Geometrically, each PDE system can be formulated as an \emph{exterior differential system} (see
\cite{BCG}).
In particular, since the current work will mainly be concerned with hyperbolic Monge-Amp\`ere systems, 
we make the definition below, following \cite{BGH1}.

\begin{Def} A \emph{hyperbolic Monge-Amp\`ere system} is an exterior differential system  $(M,\I)$, where $M$ is a 5-manifold;
	$\I\subset \Omega^*(M)$ is a differential ideal satisfying:
	for each $p\in M$, there exists an open neighborhood $U$ of $p$ on which 
	$\I|_U$ is algebraically generated by
	\[
	\theta\in\I\cap \Omega^1(U) \qquad\text{and} \qquad \ed\theta, \Omega\in
	\I\cap \Omega^2(U);
	\]
	 these generators must satisfy:
		\begin{enumerate}[\qquad (1)]
			\item{$\theta\W (\ed\theta)^2 \ne 0$;}
		
			\item{pointwise, ${\rm span}\{ \ed\theta, \Omega\}$, modulo $\theta$, has rank $2$;}
		
			\item{$(\lambda \ed\theta+ \mu \Omega)^2\equiv 0 \mod \theta$ has two distinct solutions $[\lambda_i:\mu_i]\in \R\P^1$ $(i = 1,2)$. }
		\end{enumerate}		
\label{hyperbolicMAdefinition}
\end{Def}
		
For example, the sine-Gordon equation \eqref{SG} may be formulated as a hyperbolic Monge-Amp\`ere
system $(M,\I)$ in the following way. Let $M\subset \R^5$ be an open domain with coordinates $(x,y,u,p,q)$.
Let $\I\subset \Omega^*(M)$ be the ideal generated by
			\[
				\theta = \ed u - p \ed x - q \ed y, \quad \ed\theta = \ed x \W \ed p + \ed y\W \ed q
			\]
and
			\[
				 \Omega = \(\ed p -\frac{1}{2}\sin(2u)\ed y \)\W \ed x.	
			\]
It is easy to see that a surface $\iota: S\hookrightarrow M$ satisfying $\iota^*(\ed x\W\ed y)\ne 0$
corresponds to a solution of 
\eqref{SG} if and only if $\iota^*\I = 0$. 
Moreover, the three conditions in Definition \ref{hyperbolicMAdefinition}
hold.

In general, let $(M,\I)$ be a hyperbolic Monge-Amp\`ere system. It is easy to show that, each point $p\in M$
	 has an open neighborhood $U$ on which 
	 there exist $1$-forms $\theta,\omega^1,\omega^2,\omega^3,\omega^4$, everywhere linearly independent, such that\footnote{Here $\<\cdots\>_\alg$ denotes the ideal 
	 in $\Omega^*(M)$ generated algebraically
by the differential forms enclosed in the brackets.
Without the subscript `\alg', $\<\cdots\>$ denotes the algebraic ideal generated by
the differential forms in the bracket and their exterior derivatives.
}
	 			\[
					\I|_U = \<\theta, \omega^1\W\omega^2, \omega^3\W\omega^4\>_{\alg}.
				\]		
Such a list of $1$-forms is called a \emph{$0$-adapted (local) coframing of $(M,\I)$}.
				
Suppose that $(\theta,\omega^1,\cdots,\omega^4)$ is a $0$-adapted local coframing of
a hyperbolic Monge-Amp\`ere system $(M,\I)$. Let $\iota: S\hookrightarrow M$
be an embedded surface satisfying $\iota^*\I = 0$ and the independence condition 
$\iota^*(\omega^1\W\omega^3)\ne 0$.
Such a surface is called an \emph{integral manifold} (or, in this case, an \emph{integral surface}) of $(M,\I)$.
On $S$, the $1$-forms $\omega^1$ and $\omega^2$ are multiples of each other; so are $\omega^3$
and $\omega^4$. The equations $\omega^1 = 0$ and $\omega^3 = 0$ each define a tangent
line field on $S$.
They are just
the classical characteristics of the hyperbolic PDE corresponding to $(M,\I)$.

\begin{Def} Given a hyperbolic Monge-Amp\`ere system\footnote{The notion of a \emph{characteristic system}
												applies to hyperbolic exterior differential systems in 
												general. See \cite{BGH1}.
												}									 
$(M,\I)$, where  $\I =  \<\theta, \omega^1\W\omega^2, \omega^3\W\omega^4\>,$ the pair of Pfaffian
	systems
	  $\I_{10} = \<\theta, \omega^1,\omega^2\>$ and $\I_{01} = \<\theta, \omega^3,\omega^4\>$,
	  which are defined up to ordering, are called the
	 \emph{characteristic systems} associated to $(M,\I)$.
\end{Def}	  

We follow \cite{BGG} to give the following definition.\footnote{ 
The reader may compare this definition, which suits the current work, 
with the more general version presented in \cite{HuBacklund1}.}

\begin{Def} 
Given two hyperbolic Monge-Amp\`ere systems $(M_i,\I_i)$ $(i = 1,2)$ with
		$
			\I_i  = \<\theta_i, \ed \theta_i, \Omega_i\>_\alg,
		$
a \emph{rank-$1$ B\"acklund transformation} relating them is a $6$-manifold
	\[
		\iota: N^6\hookrightarrow M_1\times M_2
	\]		
with the natural projections $\pi_i: N\rightarrow M_i$	
satisfying:
	\begin{enumerate}[\qquad (1)]
	\item{Both $\pi_1$ and $\pi_2$ are submersions, and $\ker(\ed\pi_1)\cap \ker(\ed \pi_2) = 0$;}
	
	\item{$\lbb\pi_1^*\ed\theta_1, \pi_2^*\ed\theta_2\rbb \equiv 
		\lbb\pi_1^*\Omega_1,\pi_2^*\Omega_2\rbb \mod\pi_1^*\theta_1,\pi_2^*\theta_2$.}
	\end{enumerate}	
\label{BacklundDef}
\end{Def}

Here, $\lbb\pi_1^*\ed\theta_1, \pi_2^*\ed\theta_2\rbb$ 
denotes the subbundle of $\Lambda^2(T^*N)$ 
generated by $\pi_1^*\ed\theta_1$ and $\pi_2^*\ed\theta_2$.
(The notation $\lbb\cdots\rbb$ will be used in a similar way below.)
Condition (1) implies that,
given any integral surface $S$ of $(M_1,\I_1)$, $\pi_1^{-1}S$
is $3$-dimensional. Condition (2) implies that, on~$N$,
	\[
		\pi_1^*\lbb\ed\theta_1, \Omega_1\rbb \equiv \pi_2^*\lbb\ed\theta_2,\Omega_2\rbb
			\mod\pi_1^*\theta_1, \pi_2^*\theta_2.
	\]
It follows that, on $\pi_1^{-1}S$, 
$\pi_2^*\I_2$ restricts to be algebraically
generated by the single $1$-form $\pi_2^*\theta_2$, so the Frobenius theorem applies.
In other words, $\pi_1^{-1}S$ is foliated by a 1-parameter family of surfaces whose projections via $\pi_2$
are integral surfaces of $(M_2, \I_2)$. The same argument works when one starts with an integral surface of $(M_2, \I_2)$.

\subsection{Obstructions to Existence}			

Fix two hyperbolic Monge-Amp\`ere systems $(M,\I)$ and $(\bar M,\bar \I)$, together 
with a choice of 0-adapted coframings $(\theta,\omega^1,\cdots,\omega^4)$ and
$(\bar\theta,\bar\omega^1,\cdots,\bar\omega^4)$ defined on domains
$V\subset M$ and $\bar V\subset\bar M$, respectively. In principle, the problem
of finding a rank-$1$ B\"acklund transformation relating solutions of $(M,\I)$ with those of $(\bar M,\bar \I)$
is a problem of integration. In fact, suppose that $N^6\subset V\times \bar V$ is a 
rank-$1$ B\"acklund transformation. The condition (2) in Definition \ref{BacklundDef} implies that,
on $N$,
	\[
		\lbb\omega^1 \W \omega^2, ~\omega^3\W\omega^4\rbb\equiv 
			\lbb\bar\omega^1\W\bar\omega^2, ~\bar\omega^3\W\bar\omega^4\rbb \mod \theta,\bar\theta,
	\]
where pull-back symbols are dropped for clarity.	

One can always switch the pairs $(\omega^1,\omega^2)$ and $(\omega^3,\omega^4)$, 
if needed, to arrange that
	\[	\begin{split}
		\lbb\omega^1\W\omega^2\rbb&\equiv \lbb\bar\omega^1\W\bar\omega^2\rbb \mod \theta,\bar\theta,
\\
		\lbb\omega^3\W\omega^4\rbb&\equiv\lbb\bar\omega^3\W\bar\omega^4\rbb
		\mod \theta,\bar\theta.
		\end{split}	\]	
It follows that, on $N$, there exist 
functions $s^i, t^{i}, u^{i}_j $ $(i,j = 1,\ldots,4)$ such that
	\begin{align*}
		\bar\omega^i = s^i\theta + t^i\bar\theta + u^i_j \omega^j
	\end{align*}
with $u^\alpha_r = u^s_\beta = 0$, $\det(u^\alpha_\beta), \det(u^r_s)\ne 0$ 
$(\alpha,\beta = 1,2;r,s = 3,4)$. Conversely, the existence of a rank-$1$ B\"acklund transformation reduces to analyzing the integrability of a Pfaffian system $(P, \J)$, where
$P$ is the product 
$M\times \bar M\times U$ with $U$ being a domain for the parameters $s^i,t^i$ and $u^i_j$; $\J$ is
the differential ideal generated by the four $1$-forms
	\[
		\bar\omega^i - (s^i\theta + t^i\bar\theta + u^i_j \omega^j), \quad i = 1,\ldots,4.
	\]

Theoretically, this is a type of problem that can be handled by the Cartan-K\"ahler theory 
(see \cite{BCG}). However, direct application of the idea above seems to have limited value
for two reasons. 
One, `fixing two hyperbolic Monge-Amp\`ere systems' is not general enough;
and the choice of $0$-adapted coframings is quite arbitrary. Two, the calculation 
involved in analyzing the Pfaffian system
$(P,\J)$ often quickly becomes enormous and difficult to manage. 
	
On the other hand, we can still analyze the existence of rank-$1$ B\"acklund transformations
as an integrability problem. Here we take an analogous but different approach than the one described
above. Instead of considering a specific pair of hyperbolic Monge-Amp\`ere systems with chosen 
adapted coframings, we consider the respective $G$-structure bundles, say, $\mathcal{G}$
and $\bar{\mathcal{G}}$ associated to the Monge-Amp\`ere systems. 
We then establish a rank-$4$ Pfaffian system $\mathcal{J}'$
on $\mathcal{G}\times \bar{\mathcal{G}} \times U'$, where $U'$ is a parameter 
space. Suitable integral manifolds of $\J'$ correspond to the desired B\"acklund transformations. 

This approach, based on the $G$-structure bundles instead of the manifolds of hyperbolic Monge-Amp\`ere
systems, allows one to work with all hyperbolic Monge-Amp\`ere systems at the same time.
By not specifying a coframing, it is possible to express integrability conditions in terms of 
the invariants of the Monge-Amp\`ere systems, leading to new `obstruction-to-existence' results. 

This addresses the `generality issue' mentioned above. Yet the magnitude of the calculation remains
a challenge. Being aware of this, we have assumed, for a significant portion of this work,
that the hyperbolic Monge-Amp\`ere systems under consideration
are both Euler-Lagrange
 (Section 4), which is the case for many known examples of B\"acklund transformations.
Furthermore, at a certain stage, we assume that the rank-$1$ B\"acklund transformations are of a particular type, which
we call `special' (Section 5). Such B\"acklund transformations can be divided into 4 subtypes, which we name
as {\bf Type I, IIa, IIb} and {\bf III}. 
Under these assumptions, we obtain our main obstruction results:
\vskip 2mm
{\noindent{\bf Proposition \ref{SpecialBackLemma}.}\emph{
	$\bf(A)$ If a pair of hyperbolic Euler-Lagrange systems are related by a Type I special
	B\"acklund transformation, then one of them must be positive, the other negative.
	\\
$\bf(B)$
Two hyperbolic Euler-Lagrange systems related by a Type III special B\"acklund transformation cannot be both degenerate.}
}
\vskip 2mm
{\noindent{\bf Theorem \ref{ThmIIa}.}
	\emph{If two hyperbolic Euler-Lagrange systems are related by a Type IIa special rank-$1$ 
			B\"acklund transformation, then each of them corresponds (up to contact equivalence) 
			to a second order PDE of the form 
								$
									z_{xy} = F(x,y,z,z_x,z_y).
								$			
	}
}
\vskip 2mm
{\noindent{\bf Theorem \ref{PropIIb}.}
	\emph{Let $(M,\I)$ and $(\bar M,\bar\I)$ be two hyperbolic Euler-Lagrange systems.
					If $\phi:N\rightarrow M\times \bar M$ defines a Type IIb special rank-$1$ B\"acklund transformation relating $(M,\I)$ and $(\bar M,\bar \I)$, 
					then each of $(M,\I)$ and $(\bar M,\bar\I)$ must have 
					a characteristic system that contains a rank-1 integrable subsystem.	
	}
}
\vskip 3mm
In Section \ref{newExamples}, we discover examples of Type III special rank-$1$ B\"acklund transformations
relating a degenerate hyperbolic Euler-Lagrange system with a non-degenerate one.

We provide a list of open questions in Section \ref{openquestions}.

\section{First Monge-Amp\`ere Invariants \label{InvMA}}

Let $(M,\I)$ be a hyperbolic Monge-Amp\`ere system. 
Let $\mathcal{G}_0$ be the principal bundle over $M$ consisting of $0$-adapted 
coframes of $(M,\I)$. It is easy to verify that $\mathcal{G}_0$ is a
$G_0$-structure, where $G_0\subset \GL(5,\R)$ is a $13$-dimensional Lie subgroup.
In \cite{BGG}, the reduction of $\mathcal{G}_0$ to a $G_1$-structure $\mathcal{G}_1$ is performed
such that the following structure equations hold on $\mathcal{G}_1$:
			\begin{align}
				\ed\left(
					\begin{array}{c}
						\omega^0\\
						\omega^1\\
						\omega^2\\
						\omega^3\\
						\omega^4
					\end{array}
				\right) &= -\left(
							\begin{array}{ccccc}
								\phi_0 &0&0&0&0\\
								0&\phi_1&\phi_2&0&0\\
								0&\phi_3&\phi_4&0&0\\
								0&0&0&\phi_5&\phi_6\\
								0&0&0&\phi_7&\phi_8
							\end{array}	
						\right)\W
				\left(
					\begin{array}{c}
						\omega^0\\
						\omega^1\\
						\omega^2\\
						\omega^3\\
						\omega^4
					\end{array}
				\right)  \label{StrEqnMA}\\
				&\qquad\qquad\qquad
				+ 	\left(	
				\begin{array}{c}
				\omega^1\W\omega^2+\omega^3\W\omega^4\\
				(V_1+V_5)\omega^0\W\omega^3+(V_2+V_6)\omega^0\W\omega^4\\
				 (V_3+V_7)\omega^0\W\omega^3+(V_4+V_8)\omega^0\W\omega^4\\
				 (V_8-V_4)\omega^0\W\omega^1+(V_2-V_6)\omega^0\W \omega^2\\
				 (V_3-V_7)\omega^0\W\omega^1+(V_5-V_1)\omega^0\W\omega^2
				\end{array}
				\right),\nonumber
			\end{align}
		where $\phi_0  = \phi_1+\phi_4 = \phi_5 + \phi_8$, and $G_1\subset G_0$ is the subgroup generated
		by
		\begin{equation}
			g = \left(\begin{array}{ccc}a&{\bf 0}&{\bf 0}\\
							{\bf 0}&A&0\\
							{\bf 0}&0&B
					\end{array}\right), ~A,B\in \GL(2,\R), ~a = \det(A) = \det(B), \label{1adaptedg}
		\end{equation}
		and
		\begin{equation}			
		 J = \left(\begin{array}{ccc}1&{\bf 0}&{\bf 0}\\
							{\bf 0}&0&I_2\\
							{\bf 0}&I_2&0
				\end{array}\right)\in \GL(5,\R).	\label{defmatrixJ}
		\end{equation}  

\begin{Def} Let $(M,\I)$ be a hyperbolic Monge-Amp\`ere system.  A
	 \emph{$1$-adapted coframing}
	 of  $(M,\I)$ with domain $U\subset M$ 
	is a section $\bs{\eta}: U\rightarrow \mathcal{G}_1$.
\end{Def}

Following \cite{BGG}, we introduce the notation\footnote{To be precise, these $S_i$ are $1/2$ times those defined
in \cite{BGG} with the same notation.} 
\begin{equation}
	S_1:=\left(\begin{array}{cc} V_1&V_2\\V_3&V_4\end{array}\right), \quad
	S_2:=\left(\begin{array}{cc} V_5&V_6\\V_7&V_8\end{array}\right).		\label{S1S2def}
\end{equation}
It is shown in \cite{BGG} that
\begin{Prop} \label{VTrans}
Along each fiber of $\mathcal{G}_1$, 
			\begin{equation}
				S_i(u\cdot g) = a A^{-1} S_i(u) B, \quad (i = 1,2)	\label{S1S2trans}
			\end{equation}
for any $g = \diag(a;A;B)$ in the identity component of $G_1$. Moreover,
			\begin{equation}
				S_1(u\cdot J) = \left(\begin{array}{cc} -V_4&V_2\\V_3&-V_1\end{array}\right), \quad S_2(u\cdot J) = \left(\begin{array}{cc} V_8&-V_6\\-V_7&V_5\end{array}\right).	\label{S1S2trans2}
			\end{equation} 
\end{Prop}				

Proposition \ref{VTrans} has a simple interpretation: the matrices $S_1$ and $S_2$
		correspond to two invariant tensors under the $G_1$-action. More explicitly,
		one can verify that the quadratic form
		\begin{equation}
			\Sigma_1:= V_3~\omega^1\omega^3-V_1~\omega^1\omega^4+V_4~\omega^2\omega^3 - 
						V_2~\omega^2\omega^4
		\end{equation}
		and the $2$-form
		\begin{equation}
			\Sigma_2:=V_7~\omega^1\W\omega^3 - V_5~\omega^2\W\omega^3+V_8~\omega^1\W\omega^4
						 - V_6~\omega^2\W\omega^4
		\end{equation}
		are $G_1$-invariant, which implies that $\Sigma_1, \Sigma_2$ 
		are locally well-defined on $M$.

An infinitesimal version of Proposition \ref{VTrans} will be useful: for $i = 1,2$,
			\begin{equation}
				\ed S_i \equiv \left(\begin{array}{cc} \phi_4 &-\phi_2\\-\phi_3 &\phi_1\end{array}\right) S_i +
							 S_i \left(\begin{array}{cc} \phi_5&\phi_6\\\phi_7&\phi_8\end{array}\right) \mod \omega^0,\omega^1,\ldots,\omega^4.
				\label{VonFiber}			 
			\end{equation}

An \emph{Euler-Lagrange system}, in the classical calculus of variations, is a system of PDEs whose solutions
correspond to the stationary points of a given first-order functional. In \cite{BGG}, an Euler-Lagrange system
is formulated as a Monge-Amp\`ere system\footnote{See Definitions 1.3 and 1.4 of \cite{BGG}}; 
moreover, it is shown:

\begin{Prop} {\rm (\cite{BGG})} A hyperbolic Monge-Amp\`ere system is locally equivalent to an Euler-Lagrange system if and only if $S_2$ vanishes. \label{EL}\end{Prop}

{\remark Proposition \ref{EL} says that the property of 
being \emph{Euler-Lagrange} is intrinsic. From now on, we will treat this Proposition as our
`definition' of hyperbolic Euler-Lagrange (Monge-Amp\`ere) systems.}

\begin{Prop} {\rm (\cite{BGG})} A hyperbolic Monge-Amp\`ere system corresponds to the wave equation $z_{xy}=0$ (up to contact equivalence) if and only if $S_1 = S_2 = \bs 0$. \label{wave}\end{Prop}

The following result will also be useful.

\begin{Prop} A hyperbolic Monge-Amp\`ere system $(M,\I)$, where $\I$ is algebraically generated by 
$\theta$, $\omega^1\W\omega^2$ and $\omega^3\W\omega^4$, locally corresponds to 
a PDE of the form $z_{xy} = F(x,y,z,z_x,z_y)$ (up to contact equivalence) 
if and only if each of the characteristic systems 
$\I_{10} =\< \theta,\omega^1,\omega^2\>$ and $\I_{01} = \<\theta,\omega^3,\omega^4\>$
admits a rank-$1$ integrable subsystem. \label{zxyProp}
\end{Prop}

\emph{Proof.} One direction is immediate. In fact, formulating the PDE $z_{xy} = F(x,y,z,z_x, z_y)$ 
as a hyperbolic Monge-Amp\`ere system, one easily notices that $\ed x$ and $\ed y$, respectively belonging
to the two characteristic systems,
 are integrable. 
 
  For the other direction, assume that $(M,\I)$ has the property that each of
 $\I_{10}$ and $\I_{01}$ has a rank-$1$ integrable subsystem; and let $(\theta,\omega^1,\ldots,\omega^4)$
 be a coframing defined on a domain $U\subset M$ satisfying 
 \[
 	\ed \theta \equiv \omega^1\W\omega^2+\omega^3\W\omega^4 \mod \theta.
\] 
 For $\I_{10}$, this means that a certain
 linear combination $A\theta+B\omega^1+C\omega^2$, where $A,B,C$ (not all zero) are functions on $U$,
 is closed; hence, by shrinking $U$, if needed, 
we can find a function $x$ defined on $U$ so that $\ed x = A\theta+B\omega^1+C\omega^2$. 
 Since $\theta$ is a contact form, $B,C$ cannot both be zero. Without loss of generality,
 assume that $B\ne 0$.  Let $\hat\omega^1 = A\theta+B\omega^1+C\omega^2 = \ed x$
 and $\hat \omega^2 = (1/B)\omega^2$. Similarly, there 
 exist functions $A',B',C'$ (assuming $B'\ne 0$) and $y$
 such that $\hat\omega^3 = A'\theta+B'\omega^3+C'\omega^4 = \ed y$. Let $\hat\omega^4 = 
 (1/{B'})\omega^4$.
 
 Now we have 
 \[
 	\ed\theta \equiv \ed x\W \hat\omega^2 + \ed y\W\hat\omega^4\mod \theta;
\]	 
hence, the system
 $\<\theta,\ed x,\ed y\>$ is completely integrable. By the Frobenius theorem, 
 there exists a function $z$ such that $\<\theta,\ed x,\ed y\> = \<dz, \ed x,\ed y\>$.
 In other words, there exist functions $g, p,q$ $(g\ne 0)$ defined on $U$ such that 
 \[
 	(1/g) \theta= \ed z - p\ed x - q\ed y.
\]  
 This implies that 
 \[
 	\ed x\W \ed p + \ed y\W \ed q \equiv  (1/g)(\ed x\W\hat\omega^2+\ed y\W\hat\omega^4) \mod\theta.
\]
 By Cartan's Lemma, there exists a function $F$ such that 
 	\[
		\hat\omega^2 \equiv g\ed p - gF \ed y \mod \ed x,\theta; \qquad \hat\omega^4\equiv g\ed q-gF \ed x \mod \ed y,\theta.
	\]
The vanishing of $\theta$ and $\hat\omega^1\W\hat\omega^2$ on integral surfaces then implies that locally
the corresponding Monge-Amp\`ere equation is equivalent to $z_{xy} = F(x,y,z,z_x,z_y)$.	\qed
\vskip 2mm 
Now we turn to hyperbolic Euler-Lagrange systems. 

By Proposition \ref{VTrans}, the sign of $\det(S_1)$ is independent of the choice of 1-adapted coframings. Hence, each hyperbolic Euler-Lagrange system belongs to exactly one of the following three classes.

\begin{Def} Given a hyperbolic Euler-Lagrange system $(M,\I)$, it is said to be
	\begin{itemize}
		\item{\emph{positive} if $\det(S_1)>0$;}
	
		\item{\emph{negative} if $\det(S_1)<0$;}
	
		\item{\emph{degenerate} if $\det(S_1) = 0$.}
	\end{itemize}	
\end{Def}

{\Example \label{K=-1ex}

The oriented orthonormal frame bundle $\mathcal{O}$ over the Euclidean space $\E^3$ consists
of elements of the form $(\x,\e_1,\e_2,\e_3)$, where $\x\in \E^3$, and $(\e_1, \e_2, \e_3)$
is an oriented orthonormal frame at $\x$. On $\mathcal{O}$, define the $1$-forms $\omega^i$ and $\omega^i_j$ by
		\[
			\ed\x = \e_i\omega^i,\quad \ed \e_j = \e_i\omega^i_j.
		\]
We have the standard structure equations:
			\begin{align*}
				\ed \omega^i & = -\omega^i_j \W \omega^j,\\
				\ed \omega^i_j & = -\omega^i_k \W \omega^k_j,
			\end{align*}
where $\omega^i_j+\omega^j_i = 0$. Consider the natural quotient 
		\begin{align*}
			\pi: \mathcal{O}&\rightarrow M:= \E^3\times\S^2\\
			 (\x,\e_1,\e_2,\e_3)&\mapsto (\x,\e_3).
		\end{align*}
The differential forms $\omega^3, \ed\omega^3, \ed \omega^1_2 + \omega^1\W\omega^2$
are annihilated by $\ker(\ed \pi)$ and are invariant along the fibres of $\pi$. Therefore,
they are defined on $M$. The exterior differential system $(M,\J)$, where
		\[
			\J = \<\omega^3, \ed\omega^3, \ed \omega^1_2 + \omega^1\W\omega^2\>_\alg,	
		\]
is hyperbolic Monge-Amp\`ere. Its integral surfaces correspond to generalized surfaces in $\E^3$ with Gauss 
curvature $K = -1$. 	 
		
Now consider the change of basis
			\[
				\omega^1 = -\eta^1+\eta^3,\quad
				\omega^2 = \eta^2+\eta^4,\quad
				\omega^3 = 2\eta^0,
			\]
			\[	
				\omega^1_3 = \eta^2 - \eta^4,\quad
				\omega^2_3 = \eta^1+\eta^3. 
			\]
In terms of the $\eta^i$ $(i = 0,\ldots,4)$, we have $\J=\<\eta^0,\eta^1\W\eta^2,\eta^3\W\eta^4\>_\alg$. 
It can be verified that $\bs{\eta}:= (\eta^0,\eta^1,\ldots,\eta^4)$ is 
a local $1$-adapted coframing of $(M,\J)$. Calculating using this coframing,
we obtain $V_2 = V_3 = 1$ with all other $V_i$ being zero.
It follows that $(M,\J)$ is a hyperbolic Euler-Lagrange system of the negative type.
(For a presentation using coordinates, see Appendix \ref{CalcNote}.)
}

{\Example Consider a PDE of the form $z_{xy} = f(z)$.
 (This is called an \emph{$f$-Gordon equation}.) It corresponds to a hyperbolic Euler-Lagrange system of the degenerate type, for it is easy to verify that 
			\[
				\eta^0  = \ed z - p\ed x - q\ed y,\quad
				\eta^1 = \ed x,\quad
				\eta^2 = \ed p - f(z) \ed y,
			\]
			\[	
				\eta^3 = \ed y,\quad
				\eta^4 = \ed q - f(z) \ed x
			\]
	form a 1-adapted coframing of the corresponding hyperbolic Monge-Amp\`ere system. 
	Using this coframing, one can calculate that $V_3 = -f'(z)$ with all other $V_i$ being identically zero.	}

{\Example \label{K=1ex}
Consider the Lorentzian space $\E^{2,1}$. Let $\pi:\mathcal{F}\rightarrow \E^{2,1}$
be the oriented pseudo-orthonormal frame bundle (with $\SO(2,1)$ fibres) 
consisting of $(\x,\e_1,\e_2,\e_3)$ satisfying
	\[\begin{split}
		&\x\in \E^{2,1}, \quad \e_i\in \R^3, \\
		 \e_i\cdot\e_j = 0~ (i\ne j),  
		&\quad \e_1\cdot\e_1 = -\e_2\cdot\e_2= \e_3\cdot\e_3 = 1,
	\end{split}	
	\]
where $\cdot$ stands for the inner product on $\E^{2,1}$ with the signature $(+,-,+)$.	

Define the $1$-forms $\omega^i, \omega^i_j$ on $\mathcal{F}$ by
		\[
			\ed \x = \e_i\omega^i,\qquad \ed \e_j = \e_i\omega^i_j.
		\]
We have the structure equations
		\begin{align*}
				\ed \omega^i & = -\omega^i_j \W \omega^j,\\
				\ed \omega^i_j & = -\omega^i_k \W \omega^k_j,
		\end{align*}
with $\omega^2_\alpha = \omega^\alpha_2$~ $(\alpha = 1, 3)$ and $\omega^1_3 +\omega^3_1 = 0$.
Consider the quotient
		\begin{align*}
			\pi: \mathcal{F}&\rightarrow M:= \E^{2,1}\times \mathcal{H}\\
			 (\x,\e_1,\e_2,\e_3)&\mapsto (\x,\e_3),
		\end{align*}
where $\mathcal{H}$ is the hyperboloid in $\R^3$ defined by $x^2 -y^2 + z^2 = 1$.

The differential forms $\omega^3, \ed\omega^3,\omega^1\W\omega^2 -\omega^1_3\W\omega^2_3$
are annihilated by $\ker(\ed \pi)$ and are invariant along the fibres of $\pi$. Thus they are defined
on $M$. The exterior differential system $(M,\J)$, where
	\[
		\J = \<\omega^3,\ed\omega^3,\omega^1\W\omega^2 -\omega^1_3\W\omega^2_3\>_\alg
	\]
is hyperbolic Monge-Amp\`ere. Its integral surfaces correspond to time-like (since $\e_3$, being space-like, is the normal)
surfaces in $\E^{2,1}$ with the constant Gauss curvature $K = 1$.	
	
Under the change of basis
		\[
				\omega^1 = -\eta^2+\eta^4,\quad
				\omega^2 = -\eta^1-\eta^3,\quad
				\omega^3 = -2\eta^0,
			\]
			\[	
				\omega^1_3 = -\eta^3 + \eta^1,\quad
				\omega^2_3 = \eta^4+\eta^2,
			\]
one can verify that $\bs\eta:=(\eta^0,\eta^1,\ldots,\eta^4)$ is a local $1$-adapted coframing of
$(M,\J)$. Computing using this coframing, we find that $V_2 = 1$, $V_3 = -1$, all other $V_i$ being zero.
It follows that $(M,\J)$ is a hyperbolic Euler-Lagrange system of the positive type.
(For a presentation using coordinates, see Appendix \ref{CalcNote}.)

}

{\remark
One can verify that the hyperbolic Monge-Amp\`ere systems occurring in Clelland's classification
\cite{C01} of homogeneous rank-$1$ B\"acklund transformations are Euler-Lagrange.
}

\begin{Prop} A hyperbolic Euler-Lagrange system
of either the positive or the negative type
is \emph{not} contact equivalent to any PDE of the form $z_{xy} = F(x,y,z,z_x,z_y)$.
\end{Prop}

\emph{Proof.} Given a hyperbolic Monge-Amp\`ere system $(M,\I)$,
if $S_2 = 0$ and $S_1$
nonsingular, we can compute the \emph{derived systems}\footnote{Given a Pfaffian system
$I\subset T^*M$, its derived systems are defined
recursively: $I^{(0)} = I$, $I^{(j+1)} = \lbb\omega\in  \Gamma(I^{(j)})|~\ed\omega\equiv 0\mod
I^{(j)}\rbb\subset T^*M$.
There always exists a $k\ge 0$ such that $I^{(\ell)} = I^{(k)} =: I^{(\infty)}$ for all $\ell \ge k$. 
Any integrable subsystem of $I$ is contained in $I^{(\infty)}$.} 
(see \cite{BCG}) of $I_{10} = \lbb\omega^0,\omega^1,\omega^2\rbb$
using \eqref{StrEqnMA} and find
	\[
		I_{10}^{(1)} = \lbb\omega^1, \omega^2\rbb, \qquad I_{10}^{(2)} = 0.
	\]
It follows that $I_{10}$ has no integrable subsystem. A similar argument works for $I_{01}$. 
By Proposition \ref{zxyProp}, $(M,\I)$ cannot be contact equivalent to a PDE of the form
$z_{xy} = F(x,y,z,z_x,z_y)$.\qed

 \vskip 2mm

\section{The B\"acklund-Pfaff System \label{BacklundPfaffSection}}

In this section, we prove that, given two hyperbolic Monge-Amp\`ere systems,
the existence of a \emph{rank-$1$ B\"acklund transformation}
relating them can be formulated
as the integrability of a rank-$4$ Pfaffian system. 

We start by fixing some notations.

Let $(M,\I)$ and $(\bar M,\bar\I)$ be two hyperbolic Monge-Amp\`ere systems. 
Let ${\mathcal{G}}_1$ and $\bar{\mathcal{G}}_1$ be the respective $G_1$-structures (see Section \ref{InvMA}). 
Let $\alpha = (\alpha_0,\alpha_1,\ldots,\alpha_4)$ and
$\beta = (\beta_0,\beta_1,\ldots,\beta_4)$ be the tautological 1-forms on $\mathcal{G}_1$ and 
$\bar{\mathcal{G}}_1$, respectively. And let $\mathcal{P}$ be the product 
of $\mathcal{G}_1\times \bar{\mathcal{G}_1}$ with a space of parameters:
	\begin{align}
		\mathcal{P} &:=\left\{(u,\bar u, s,t, \mu,\epsilon)\left| ~\begin{array}{c}u\in \mathcal{G}_1,~ \bar u\in\bar{\mathcal{G}}_1;~ s,t\in \R^4;\\ \mu\ge 1;~ \epsilon = \pm 1 ~(\epsilon\mu^2\ne 1)\end{array}\right.\right\}	\label{parameterspace}\\
		&\subset \mathcal{G}_1\times \bar{\mathcal{G}_1}\times\R^9\times\{\pm 1\}.\nonumber
	\end{align}

\begin{Prop} An embedded $6$-manifold $\phi: N\rightarrow M\times\bar M$ is a rank-$1$
 B\"acklund transformation if and only if $N$ admits a
lifting $\hat\phi: N\rightarrow \mathcal{P}$ that is an integral manifold of a rank-$4$ Pfaffian system  $\mathcal{J}$ with the \emph{independence condition}
$\alpha_1\W \alpha_2\W \alpha_3\W \alpha_4\W \alpha_0 \W \beta_0\ne 0$,
 where $\J$ is generated by
		\begin{equation}
			\begin{array}{l}
			\theta_1 := \beta_1 - s_1 \alpha_0 - t_1 \beta_0 - \mu \alpha_1,\\[0.8em]
			\theta_2 := \beta_2 - s_2 \alpha_0 - t_2 \beta_0 - \mu \alpha_2,\\[0.8em]
			\theta_3 := \beta_3 - s_3 \alpha_0 - t_3 \beta_0 - \mu^{-1} \alpha_3,\\[0.8em]
			\theta_4 := \beta_4 - s_4 \alpha_0 - t_4 \beta_0 - \epsilon \mu^{-1} \alpha_4,
			\end{array}\label{backlundpfaffthetas}
		\end{equation}
at each $(u,\bar u, s_1,\ldots,s_4, t_1,\ldots,t_4,\mu,\epsilon)\in \mathcal{P}$.	
\begin{figure*}[h!]
					\begin{center}
					\begin{tikzpicture}
						\node at(0,0){$M\times\bar M$};
						\node at (0,1.5){$\mathcal{P}$};
						\node at (-2.5,0){$N^6$};
						\node at (-1.5,-1.5){$(M,\I)$};
						\node at (1.5,-1.5){$(\bar M,\bar\I)$};
						
						\draw[->] (0,1.2) -- (0,0.2); \node at (0.3,0.8) {$\tau$};
						\draw[->] (-2.1,0.2) -- (-0.3,1.3); \node at (-1.4,1) {$\hat\phi$};
						\draw[->] (-2,0) -- (-0.8,0); \node at (-1.4,-0.25) {$\phi$};
						\draw[->] (-0.3,-0.3) -- (-1.2,-1.2); \node at (-0.9,-0.6) {$\pi$};
						\draw[->] (0.3,-0.3) -- (1.2,-1.2); \node at (0.9,-0.6) {$\bar\pi$};
					\end{tikzpicture}
					\end{center}
					\label{Backlund}
\end{figure*}	
\label{BacklundPfaffProp}
\end{Prop}

 \emph{Proof.} Let $\tau, \pi, \bar \pi$ denote the obvious projections (see the diagram above). 
By the construction of $\mathcal{G}_1$ and $\bar{\mathcal{G}}_1$, on $\mathcal{P}$, we have
				\begin{align*}
					(\pi\circ\tau)^*(\I) &= \<\alpha_0,\alpha_1\W \alpha_2, \alpha_3\W \alpha_4\>,\\
					\lbb(\pi\circ\tau)^*(\theta)\rbb& \equiv \lbb\alpha_0\rbb,
				\end{align*}
and similarly for $\bar\I$	and the $\beta_i$.			
 Now assume that $\phi$ admits a lifting $\hat\phi$ that integrates the Pfaffian system $\J$. It is easy to see that, on $N$, 
 			\begin{align*}
				\hat\phi^*(\lbb\alpha_1\W \alpha_2, \alpha_3\W \alpha_4\rbb) 
				&\equiv \hat\phi^*(\lbb\beta_1\W \beta_2, \beta_3\W \beta_4\rbb) \\
										           	&\equiv \hat\phi^*(\lbb \ed\alpha_0, \ed\beta_0\rbb) \mod \hat\phi^*\alpha_0, \hat\phi^*\beta_0.
			\end{align*}
In the last congruence, we have used the assumption that $\mu^2\ne \epsilon$, which guarantees that 
the bundle $\lbb\alpha_1\W \alpha_2+\alpha_3\W \alpha_4, ~\beta_1\W \beta_2+\beta_3\W \beta_4\rbb$ has rank $2$ modulo $\alpha_0,\beta_0$ when pulled back to $N$.	
It follows that $\phi: N\rightarrow M\times\bar M$ defines a rank-$1$ B\"acklund transformation.		

Conversely, suppose that $\phi: N\rightarrow M\times \bar M$ defines a rank-$1$ B\"acklund transformation. Let $\bs{\eta}:=(\eta_0,\eta_1,\ldots,\eta_4)$ (resp. $\bs{\xi}:=(\xi_0,\xi_1,\ldots,\xi_4)$) be a local 1-adapted coframing defined
on a domain in $M$ (resp. $\bar M$).
We have, by definition,
			\[
			 	\phi^*(\lbb\eta_1\W\eta_2, \eta_3\W\eta_4\rbb) \equiv \phi^*(\lbb\xi_1\W \xi_2, \xi_3\W\xi_4\rbb) \mod \phi^*\eta_0, \phi^*\xi_0.
			\]
By switching the pairs $(\eta_1,\eta_2)$ and $(\eta_3,\eta_4)$ and by shrinking $N$, if needed, we can assume that, on $N$,
			\begin{align*}
				\lbb\eta_1\W\eta_2 \rbb&\equiv \lbb\xi_1\W \xi_2\rbb \mod \eta_0,\xi_0,\\
				\lbb\eta_3\W\eta_4\rbb& \equiv \lbb\xi_3\W\xi_4\rbb  \mod\eta_0,\xi_0,
			\end{align*}
where the pull-back symbol is dropped for convenience. Consequently, there exist 16 functions $s_1,\ldots,s_4,t_1,\ldots,t_4,u_1,\ldots,u_8$ defined on $N$ such that, when restricted to $N$,
			\begin{equation}
				\begin{array}{l}
				 \xi_1 = s_1 \eta_0 + t_1 \xi_0 + u_1 \eta_1 + u_2 \eta_2,\\[0.8em]
				 \xi_2 = s_2 \eta_0 + t_2 \xi_0 + u_3 \eta_1 + u_4 \eta_2,\\[0.8em]
				 \xi_3 = s_3 \eta_0 + t_3 \xi_0 + u_5 \eta_3 + u_6 \eta_4,\\[0.8em]
				 \xi_4 = s_4 \eta_0 + t_4 \eta_0 + u_7 \eta_3 + u_8 \eta_4.
				 \end{array}	\label{xieta}
			\end{equation}
Here $u_1u_4-u_2u_3$ and $u_5 u_8-u_7u_6$ are both nonvanishing. Moreover, since
$\lbb \ed\eta_0, \ed\xi_0\rbb$ has rank $2$ modulo $\eta_0, \xi_0$, we have 
\[u_1u_4-u_2u_3\ne u_5 u_8-u_6u_7.\] 

Using the flexibility of choosing the $1$-adapted coframings $\bs\eta$ and $\bs\xi$ (see \eqref{1adaptedg}, 
\eqref{defmatrixJ}), we can
normalize some of the $u_i$. 
To be specific, we can apply $\SL(2,\R)$-actions (i.e., acting on $\bs\eta$ pointwise by a matrix $g$ in 
the form of \eqref{1adaptedg} with $a = 1$) to $(\eta_1,\eta_2)$ and $(\eta_3, \eta_4)$
to arrange that
\[
		u_2 = u_3 = u_6 = u_7 = 0
\]		
and
\[
	u_4 = \pm u_1, \quad u_8 = \pm u_5\quad (u_1,u_5>0).
\]	
Then we transform $(\eta_0,\eta_1,\ldots,\eta_4)$  pointwise by a diagonal matrix in $G_1$ to arrange that 
\[
	u_1 = u_4>0, \quad u_5 = u_1^{-1}, \quad u_8  = \epsilon u_1^{-1} \quad (\epsilon = \pm 1). 
\]
Meanwhile, this transformation scales the $s_i$ $(i = 1,\ldots,4)$. Finally, if $u_1 = u_4<1$, then we switch the pairs of indices $(1,2)$ and $(3,4)$ for $\bs\xi$ and $\bs\eta$ in $\eqref{xieta}$
and multiply the new $\eta_0,\eta_2,\eta_4$ by $\epsilon$; the resulting $u_i$ will satisfy
\[
	u_1 = u_4>1, \quad u_5 = u_1^{-1}, \quad u_8 = \epsilon u_1^{-1}.
\]

It follows that the $(\bs\eta,\bs\xi,(s_i),(t_i),u_1, \epsilon)$ thus constructed defines
a map $\hat\phi: N\rightarrow \mathcal{P}$ that is a lifting of $\phi$ and 
an integral manifold of $\J$. 
\qed
\vskip 3mm

In light of Proposition \ref{BacklundPfaffProp}, we make the following definitions.

\begin{Def} The system $(\mathcal{P},\J)$ in Proposition \ref{BacklundPfaffProp} is 
called the \emph{$0$-refined B\"acklund-Pfaff system} for rank-$1$ B\"acklund transformations relating two
hyperbolic Monge-Amp\`ere systems.
\end{Def}

\begin{Def} A $6$-dimensional
integral manifold (satisfying the 
independence condition described in Proposition \ref{BacklundPfaffProp}) of the $0$-refined B\"acklund-Pfaff system $(\mathcal{P},\J)$ is called a 
\emph{$0$-refined lifting} of the underlying B\"acklund transformation.
\end{Def}
	
In these terms, Proposition \ref{BacklundPfaffProp} says that each rank-$1$ B\"acklund transformation 
relating two hyperbolic Monge-Amp\`ere systems has a $0$-refined lifting. Of course, given such a 
B\"acklund transformation $\phi: N \rightarrow M\times\bar M$, its $0$-refined liftings
are not unique. 

\begin{lemma} Let $\phi: N\rightarrow M\times\bar M$ be a rank-$1$ B\"acklund transformation
relating two hyperbolic Monge-Amp\`ere systems. The functions $\mu\circ \hat\phi$ and $\epsilon\circ\hat \phi$ are independent of the choice of $0$-refined liftings $\hat\phi$ of $\phi$.
\label{epsilonmulemma}
\end{lemma}

\emph{Proof.} Clearly, for different choices of $\hat\phi$, the $1$-forms 
			$\phi^*\alpha_0$ and $\phi^*\beta_0$ only change by scaling. 
		On $N$, the quotient $\lambda_1/\lambda_2$ between the two solutions $\lambda_1, \lambda_2$ of the equation
					\[
						(\ed\beta_0 + \lambda \ed \alpha_0)^2 \equiv 0 \mod \alpha_0, \beta_0
					\]			
		is independent of the scaling of $\alpha_0$ and $\beta_0$ but may depend on the order of the pair
			$(\lambda_1,\lambda_2)$. 
		By \eqref{backlundpfaffthetas}, either $\lambda_1/\lambda_2$ or $\lambda_2/\lambda_1$ 
		must be equal to $\epsilon\mu^4$.
		
		When $|\lambda_1/\lambda_2| = 1$, it is necessary that $\mu = 1$ and $\epsilon = -1$.
		When $|\lambda_1/\lambda_2| \ne 1$, 
			 the ambiguity of ordering $\lambda_1$ and $\lambda_2$ can be eliminated by 
			 requiring that $|\lambda_1/\lambda_2|>1$, and hence $\lambda_1/\lambda_2
			 = \epsilon\mu^4$. 
			 It follows that $\mu$ and $\epsilon$ are both independent of the lifting.\qed
\vskip 2mm

{\remark {\bf A.}
There is a simple geometric interpretation for the two possible values of $\epsilon$. Let $(M,\I)$ and $(\bar M,\bar\I)$ be as above. Suppose that $\phi:N\rightarrow M\times \bar M$ defines a rank-$1$ B\"acklund transformation.
The 4-plane field $\mathcal{D}:=\lbb\alpha_0,\beta_0\rbb^{\perp}$ 
 on $N$ is independent of the choice of $0$-refined liftings of $\phi$. When restricted to $\mathcal{D}$, the $4$-forms
 $(\ed \alpha_0)^2$ and $(\ed \beta_0)^2$ 
define two orientations. If $\epsilon = 1$, these two orientations are the same; if $\epsilon = -1$, they are distinct.

{\bf B.}
If, for a rank-$1$ B\"acklund transformation relating two hyperbolic Monge-Amp\`ere systems, 
$\epsilon = 1$, then $\mu >1$.
This is because, if $\mu = 1$, then the condition (2) in Definition \ref{BacklundDef} will not hold.

	\label{remarkmuepsilon}}

\section{Obstructions to Integrability \label{obstructionSection}}

In this section, we express some obstructions to the integrability of $(\mathcal{P},\mathcal{J})$ 
 in terms of the invariants of the two hyperbolic Monge-Amp\`ere systems.

For convenience, we introduce new notations below.
\begin{enumerate}[\bf (A)]
	\item{Let
		$
			\eta_1:= \alpha_1, ~\eta_2:= \alpha_2,~  \eta_3:= \alpha_3, ~ \eta_4:= \alpha_4,~  \eta_5:= \alpha_0,~  \eta_6:= \beta_0.
		$}
		
		\item{The components of the pseudo-connection 1-form on $\mathcal{G}_1$ $($see \eqref{StrEqnMA}$)$ are denoted by $\varphi_0, \varphi_1,\ldots, \varphi_8$; those on $\bar{\mathcal{G}}_1$ are denoted by $\psi_0, \psi_1,\ldots, \psi_8$.}
\end{enumerate}
		
On $\mathcal{P}$, differentiating the $\theta_i$ and reducing modulo $\theta_1,\ldots,\theta_4$ yields the following congruences:
		\begin{equation}
			\ed\theta_i \equiv -\pi_{i\alpha}\W\eta_\alpha + \tau_i \quad\mod\theta_1,\ldots,\theta_4\qquad	(i = 1,\ldots,4),\label{PfaffStr}
		\end{equation}
where summation over repeated indices is intended; $\pi_{i\alpha}$ $(i = 1,\ldots,4; \alpha = 1,\ldots,6)$ are linear combinations of the 1-forms
in the set	\begin{align*}
			\mathcal{S}:= &\{\ed(s_1),\ldots, \ed(s_4); ~\ed(t_1),\ldots,\ed(t_4);~\ed\mu;\\
			& \varphi_0,\ldots,\varphi_3, \varphi_5,\varphi_6,\varphi_7;~ \psi_0,\ldots,\psi_3,\psi_5,\psi_6,\psi_7\};
		\end{align*}
the components of  the \emph{torsion} $T:=(\tau_i)$ are of the form
		\[
			\tau_i = \frac{1}{2}T_{ijk}\eta_j\W\eta_k,
		\]		
for some functions $T_{ijk}$ defined on $\mathcal{P}$ satisfying $T_{ijk}  = -T_{ikj}$.		

We use a standard method (see \cite{BCG}) to obtain from \eqref{PfaffStr} obstructions
to the integrability of $(\mathcal{P}, \mathcal{J})$. 
The key is that, on an integral manifold of $(P,\J)$, $\theta_i, \ed \theta_i$ all vanish and 
$\pi_{i\alpha}$ are linear combinations of $\eta_1,\eta_2,\ldots, \eta_6$. It follows that
an integral manifold of $(\mathcal{P}, \J)$ can only exist within the locus along which $\tau_i$ can be fully absorbed
by $\pi_{i\alpha}\W\eta_\alpha$, as long as one add suitable linear combinations of $\eta_\alpha$
to the $1$-forms in $\mathcal{S}$. 

This is precisely the approach we take in practice: we keep adjusting the $1$-forms in $\mathcal{S}$ by 
adding linear combinations of $\eta_\alpha$ to them until a point when the rank of $(\tau^i)$ cannot be further
reduced. The remaining $(\tau^i)$ must vanish along integral manifolds of $(P,\mathcal{J})$.

In our case, the matrix $(\pi_{i\alpha})$ (called the \emph{tableau}) takes the form
		\[
			(\pi_{i\alpha}) = \left(
								\begin{array}{cccccc}
									\pi_1&\pi_2&0&0&\pi_3&\pi_4\\
									\pi_5&\pi_6&0&0&\pi_7&\pi_8\\
									0&0&\pi_9&\pi_{10}&\pi_{11}&\pi_{12}\\
									0&0&\pi_{13}&\pi_{14}&\pi_{15}&\pi_{16}
								\end{array}
							\right),
		\]
where $\pi_1,\ldots,\pi_{16}$ are 1-forms\footnote{One can think of $(\pi_{i\alpha})$ as an $\mathcal{A}$-valued $1$-form, where $\mathcal{A}\subset
\R^4\otimes(\R^6)^*$ is a vector subspace. Here, the $\R$-valued $1$-forms $\pi_1,\ldots,\pi_{16}$ are just a set of 
`form coordinates' for $(\pi_{i\alpha})$.}
 linearly independent of $\theta_1,\ldots,\theta_4, \eta_1,\ldots,\eta_6$ and  among themselves.
It is easy to see that all terms in $\tau_i$ can be absorbed except the $\eta_3\W\eta_4$ terms in $\tau_1$ and $\tau_2$ and the $\eta_1\W\eta_2$ terms in $\tau_3$ and $\tau_4$.

Calculation yields
		\begin{align*}
			\ed\theta_1&\equiv  -\frac{\mu^2s_1+\epsilon t_1}{\mu^2}\eta_3\W\eta_4 \quad\mod\theta_1,\ldots,\theta_4,\eta_1,\eta_2,\eta_5,\eta_6,\\  
			\ed\theta_2&\equiv  -\frac{\mu^2s_2+\epsilon t_2}{\mu^2}\eta_3\W\eta_4 \quad\mod\theta_1,\ldots,\theta_4,\eta_1,\eta_2,\eta_5,\eta_6,
		\end{align*}
		\begin{align*}	
			\ed\theta_3&\equiv -(\mu^2 t_3-s_3)\eta_1\W\eta_2 \quad\mod\theta_1,\ldots,\theta_4,\eta_3,\eta_4,\eta_5,\eta_6,\\
			\ed\theta_4&\equiv  -(\mu^2 t_4 - s_4)\eta_1\W\eta_2 \quad\mod\theta_1,\ldots,\theta_4,\eta_3,\eta_4,\eta_5,\eta_6.
		\end{align*}
As a result, we have proved:
\begin{lemma} Integral manifolds of $(\mathcal{P},\J)$ must be contained in the locus $\mathcal{L}\subset\mathcal{P}$ defined by the equations
			\begin{align*}
				&s_3 = -t_3\mu^2,&& s_4 = -t_4\mu^2,\\
				&t_1 = -\epsilon s_1\mu^2, && t_2 = -\epsilon s_2\mu^2.
			\end{align*}
\end{lemma}	
\begin{Def} Let $\mathcal{P}_1\subset \mathcal{P}$ be the locus defined by the equations
			\[
				s_1 = t_1 = s_3 = t_3 = 0, \quad   s_4 = -t_4\mu^2, \quad  t_2 = -\epsilon s_2\mu^2.
			\]
	Let $\J_1$ be the pull-back of $\J$ to $\mathcal{P}_1$. The rank-4 Pfaffian system $(\mathcal{P}_1, \J_1)$ is called the \emph{$1$-refined B\"acklund-Pfaff system} for rank-$1$ 
	B\"acklund transformations relating two hyperbolic
	Monge-Amp\`ere systems.
\end{Def}	
\begin{Def} A $6$-dimensional integral manifold of the $1$-refined B\"acklund-Pfaff system $(\mathcal{P}_1,\J_1)$ is called a \emph{$1$-refined lifting} of the underlying B\"acklund transformation.
\end{Def}	
		
\begin{Prop} Let $(M,\I)$ and $(\bar M,\bar \I)$ be as above. 
	    Any rank-$1$ B\"acklund transformation 
	    $\phi: N\rightarrow M\times\bar M$ admits a $1$-refined lifting.
\end{Prop}
\emph{Proof.} By the previous discussion, there exists a $0$-refined lifting $\hat \phi$ of $\phi$ such that, when pulled-back via $\hat \phi$ to $N$,
				\begin{align*}
					\beta_1& = s_1 \alpha_0 - \epsilon s_1 \mu^2 \beta_0 + \mu \alpha_1,\\
					\beta_2& = s_2 \alpha_0 - \epsilon s_2 \mu^2 \beta_0 + \mu \alpha_2,\\
					\beta_3& = -t_3\mu^2 \alpha_0 + t_3 \beta_0 +\mu^{-1} \alpha_3,\\
					\beta_4& = -t_4\mu^2 \alpha_0 + t_4 \beta_0 + \epsilon \mu^{-1} \alpha_4.
				\end{align*}		
			If $s_2\ne 0$, we transform 
			the pairs $(\alpha_1,\alpha_2)^T$ and $(\beta_1,\beta_2)^T$ simultaneously 
			by
				\[
					h_1 = \left(\begin{array}{cc}  1& -{s_1}/{s_2} \\0&1\end{array}\right).
				\]	
			 If $s_2 = 0$ but $s_1\ne 0$, we
			 transform the pairs $(\alpha_1,\alpha_2)^T$ and $(\beta_1,\beta_2)^T$ simultaneously 
			by
				\[
					h_2 =  \left(\begin{array}{cc}  0& -1 \\ 1&0\end{array}\right)
				\]	
			such that the previous case applies. 
			These transformations correspond to choosing new $0$-refined liftings, and the result 
			is a $0$-refined lifting with $s_1 = 0$. 
			In a similar way, we can choose $0$-refined liftings such that, in addition, $t_3 = 0$.
			This completes the proof.
			\qed
\vskip 2mm
			
Now $\J_1$ is generated by
		\begin{equation}\begin{split}
			\theta_1 &= \beta_1  - \mu \alpha_1,\\
			\theta_2 &= \beta_2 - s_2 \alpha_0 + \epsilon s_2\mu^2 \beta_0 - \mu \alpha_2,\\
			\theta_3 &= \beta_3  - \mu^{-1} \alpha_3,\\
			\theta_4 &= \beta_4 + t_4\mu^2 \alpha_0 - t_4 \beta_0 - \epsilon \mu^{-1} \alpha_4.
			\end{split}	\label{1refinedeq}
		\end{equation}
Given a rank-$1$ B\"acklund transformation $\phi: N\rightarrow M\times\bar M$, it is easy to see that whether the product $s_2t_4$ locally vanishes is independent of the choice of $1$-refined liftings of $\phi$.
	It turns out that the case when $s_2t_4 \equiv 0$ is quite restrictive 
	when both $(M,\I)$ and
	$(\bar M,\bar\I)$ are Euler-Lagrange systems.
	
\begin{Prop} \label{s2t4=0} Let $(M,\I)$ and $(\bar M,\bar\I)$ be two hyperbolic Euler-Lagrange systems. If there exists a
rank-$1$ B\"acklund transformation $\phi: N\rightarrow M\times \bar M$ such that $s_2t_4 = 0$ 
			 on a $1$-refined lifting of $\phi$, then both $(M,\I)$ and $(\bar M,\bar\I)$ must be 
			 contact equivalent to the system representing the wave equation \[z_{xy} = 0.\]	
\end{Prop}			 
\emph{Proof.} By the Euler-Lagrange assumption and Proposition \ref{EL}, we have $S_2 = \bar S_2 = \bs 0$. Let 
					\[
							S_1 = \left(\begin{array}{cc}
										V_1&V_2\\
										V_3 &V_4
									\end{array}
									\right),\qquad 
							\bar S_1 = \left(\begin{array}{cc}
										W_1&W_2\\
										W_3 &W_4
									\end{array}
									\right).		
					\]	
			By Proposition \ref{wave}, it suffices to show that, on any integral manifold of the $1$-refined B\"acklund-Pfaff system $(\mathcal{P}_1,\mathcal{J}_1)$, we have $S_1 = \bar S_1 = \bs 0$. 
		
		   	First, we assume that $s_2 = 0$. Restricted to the locus defined by $s_2 = 0$ in $\mathcal{P}_1$, the tableau $(\pi_{i\alpha})$ associated to $(\mathcal{P}_1,\mathcal{J}_1)$ satisfies 
					\[	
							\pi_{i\alpha} = 0, \quad ( i = 1,2; ~\alpha = 3,\ldots,6).
					\]
			As a result, for $i,j = 3,\ldots,6$ and $i\ne j$, 
			the $\eta^i\W\eta^j$ terms in $\ed\theta_1$ and $\ed\theta_2$
			cannot be absorbed, and the corresponding coefficients must vanish on 
			any integral manifold of
			$(\mathcal{P}_1,\mathcal{J}_1)$. 
			Calculating with Maple\texttrademark,  we find that
					\begin{align*}
						\ed\theta_1 \equiv \mu(V_1\eta_3+V_2\eta_4)\W\eta_5 -\frac{1}{\mu}(W_1\eta_3+\epsilon W_2\eta_4)\W\eta_6+W_2t_4\mu^2\eta_5\W\eta_6,\\
						\ed\theta_2 \equiv \mu(V_3\eta_3+V_4\eta_4)\W\eta_5 - \frac{1}{\mu}(W_3\eta_3+\epsilon W_4\eta_4)\W\eta_6+W_4t_4\mu^2\eta_5\W\eta_6,
					\end{align*}
			both congruences being reduced modulo $\theta_1,\ldots,\theta_4,\eta_1,\eta_2$. 
			It follows that $S_1 = \bar S_1 = \bs 0$.
			
			The case when $t_4=0$ is similar. \qed		
							
\begin{Prop} \label{ELTorsion}
On any $1$-refined lifting of a rank-$1$ B\"acklund transformation relating two hyperbolic Euler-Lagrange systems, the following expressions must vanish:
			\begin{align*}
				\Phi_1 & := -\mu^4V_1+\epsilon W_1,\qquad
				\Phi_2 := -\mu^4V_2 +W_2,\\
				\Phi_3 & := \mu^4 W_4 - V_4,\qquad\quad\phantom{'}
				\Phi_4 := \mu^4 W_2 - V_2.
			\end{align*}
\end{Prop}	
\emph{Proof.} This is evident when such a B\"acklund transformations satisfies $s_2t_4 = 0$. In fact, by
 Proposition \ref{s2t4=0}, the functions $V_i$ and $W_i$ vanish identically on $\mathcal{P}_1$.
			Otherwise, restricting to the open subset of $\mathcal{P}_1$ defined by $s_2t_4\ne 0$, a calculation similar to that in the proof of Proposition \ref{s2t4=0} shows that the torsion of the Pfaffian system $\J_1$ can be absorbed only when $\Phi_i$ $(i = 1,\ldots,4)$ 
			are all zero. \qed

\begin{corollary} If two hyperbolic Euler-Lagrange systems are related by a rank-$1$ B\"acklund transformation with $\epsilon = 1$, then they are either both degenerate or both nondegenerate.	\label{CorEps1}\end{corollary}
\emph{Proof.} We have noted above (see Remark \ref{remarkmuepsilon}) that, if a rank-$1$ B\"acklund transformation relating two hyperbolic Monge-Amp\`ere systems
 satisfies $\epsilon = 1$, then $\mu > 1$. The vanishing of $\Phi_2$ and $\Phi_4$ on a $1$-refined lifting of such a B\"acklund transformation then implies that, 
				on such a lifting,
				\[
					V_2 = W_2 = 0.
				\]
By the vanishing of $\Phi_1$ and $\Phi_3$, it is easy to see
that the matrices $S_1$ and $\bar S_1$ are either both degenerate or both nondegenerate. 	\qed
\\

Note that, in the proof of Corollary \ref{CorEps1}, the condition $V_2 = W_2 = 0$ is meaningful only if it is independent of the choice of
$1$-refined liftings of a B\"acklund transformation. To make this point explicit, we prove
the following proposition, which 
shows how $1$-refined liftings of a rank-$1$ B\"acklund transformation relate to each other.

\begin{Prop}\label{lev1Perm}
 Let $\phi$ be a rank-$1$ B\"acklund transformation relating two hyperbolic Monge-Amp\`ere systems satisfying $s_2t_4\ne 0$ on its $1$-refined liftings. 
There exists a subgroup $H\subset G_1\times {G}_1$
such that any two $1$-refined liftings of $\phi$ are related, in the $\mathcal{G}_1\times\bar{\mathcal{G}}_1$ component, by an ${H}$-valued transformation. Moreover,
	\begin{enumerate}[\bf i.]
	\item{If $\mu > 1$, then $H$ is the subgroup generated by elements of the form $h = (g,g')$ where
						\begin{equation}
								g = 
						\left(
							\begin{array}{ccccc}
								h_0 &0&0&0&0\\
								0&h_1&0&0&0\\
								0&h_3&h_0h_1^{-1}&0&0\\
								0&0&0&h_2&0\\
								0&0&0&h_4&h_0h_2^{-1}
							\end{array}	
						\right)\in G_1,		\label{formofg}
						\end{equation}
			and $g'$ is the result of replacing $h_4$ in $g$ by $\epsilon h_4$.		}	
	\item{If $\mu = 1$, then $H$ is generated by the subgroup in case {\bf i} and the element
						\[
							h_J = (J',J), 
						\]
			where $J$ is as in \eqref{defmatrixJ} and $J' = \diag(-1,1,-1,1,-1)$J.}
	\end{enumerate}					
\end{Prop}

\emph{Proof.} By \eqref{1refinedeq}, any two $1$-refined liftings of $\phi$ must be related by a
		pointwise $H$-action, where
		$H\subset G_1\times G_1$ is a subgroup.
		
		Note that $G_1$ has two connected components: $G_e$ consisting of those elements
		of the form \eqref{1adaptedg} and $G_1\backslash G_e = JG_e$. 
		
		When $\mu>1$, one must have $H\subset G_e\times G_e$. This is because switching 
		$\lbb\alpha_1,\alpha_2\rbb$ with $\lbb\alpha_3,\alpha_4\rbb$
		enforces switching $\lbb\beta_1,\beta_2\rbb$ with $\lbb\beta_3,\beta_4\rbb$,
		which transforms $\mu$ into $\mu^{-1}$, which is not permissible. Now suppose 
		that $h = (g,g')\in H$.
		It is easy to see, by \eqref{1refinedeq} and the assumption $s_2t_4\ne 0$ that 
		$g\in H$ must be of the form \eqref{formofg}. The form of $g'$ is then determined, as stated in {\bf i}.
		
		When $\mu = 1$, we have $\epsilon = -1$. 
		In this case, we do allow switching between the indices $(1,2)$ and $(3,4)$, by 
		an action of $J$ on the $\beta$'s and a corresponding action on the $\alpha$'s so that
		the form of \eqref{1refinedeq} is preserved.
		
		Finally, note that a pointwise action by $H$ described above transforms a $1$-refined lifting of $\phi$
		to a $1$-refined lifting of $\phi$.
		
		This completes the proof.\qed
		
{\remark		
	By \eqref{S1S2trans} and \eqref{S1S2trans2}, 
	$h = (g,g')$ and $h_J = (J',J)$ in Proposition \ref{lev1Perm} act on $S_1 = (V_i)$ and $\bar S_1 = (W_i)$ in the following way:
	\begin{align*}
		S_1(u\cdot g) &= h_0\(\begin{array}{cc}
			h_1^{-1}&0\\
			-h_0^{-1}h_3&h_0^{-1}h_1
		\end{array}\)S_1(u)\(\begin{array}{cc}
			h_2&0\\
			h_4&h_0h_2^{-1}
		\end{array}\),\\
		\bar S_1(\bar u\cdot g') &=  h_0\(\begin{array}{cc}
			h_1^{-1}&0\\
			-h_0^{-1}h_3&h_0^{-1}h_1
		\end{array}\)\bar S_1(\bar u)\(\begin{array}{cc}
			h_2&0\\
			\epsilon h_4&h_0h_2^{-1}
		\end{array}\),\\
		S_1(u\cdot J')& = \left(\begin{array}{cc}
			V_4&V_2\\
			V_3&V_1
		\end{array}\right)(u),\quad
		\bar S_1(\bar u\cdot J) = \left(\begin{array}{cc}
			-W_4&W_2\\
			W_3&-W_1
		\end{array}\right)(\bar u).
	\end{align*}
In particular, if $W_2 =V_2 = 0$ holds on a $1$-refined lifting of $\phi$, 
then it holds on any other $1$-refined lifting of $\phi$.	\label{Hremark}
}

\section{A Special Class of B\"acklund Transformations	\label{specialclasssection}}

In the previous section, we have seen that, to a rank-$1$ B\"acklund transformation $\phi: N\rightarrow M\times \bar M$ relating two hyperbolic Monge-Amp\`ere systems $(M,\I)$ and $(\bar M,\bar\I)$,
we can associate a function $\mu: N\rightarrow [1,\infty)$ that is 
independent\footnote{In fact, it is easy to see, from the point of view of \cite{HuBacklund1} 
that $\mu$ is a local invariant
of the B\"acklund transformation.}
of the $0$-refined liftings of $\phi$.

\begin{Def}
A rank-$1$ B\"acklund transformation relating two hyperbolic Monge-Amp\`ere systems is said to be \emph{special} if $\mu = 1$.
\end{Def}

For the rest of this section, we will focus on special B\"acklund transformations relating two hyperbolic
Euler-Lagrange systems. A motivation for this is that many classical B\"acklund transformations are of this type
(cf. \cite{C01}, \cite{C13}).

By Proposition \ref{ELTorsion}, given a special rank-$1$ B\"acklund transformation $\phi:N\rightarrow M\times\bar M$ relating two hyperbolic Euler-Lagrange systems, the following equalities must hold  
on any $1$-refined lifting of $\phi$:
		\begin{equation}
			 \epsilon = -1, \quad W_1 = - V_1,\quad W_2 = V_2, \quad W_4 = V_4.   \label{level1ELrel}
		\end{equation}	
By Remark \ref{Hremark}, 
these conditions are invariant under the $H$-action defined in Proposition \ref{lev1Perm}. 
		
Now let $\mathcal{P}_s\subset\mathcal{P}_1$ be defined by the equations 
		\[
			\mu = 1, \qquad \epsilon = -1.
		\]	
				
\begin{Prop} 
Any $1$-refined lifting of a special rank-$1$ B\"acklund transformation relating two hyperbolic Euler-Lagrange systems
 is completely contained in $\mathcal{P}_s$. Moreover, on such a lifting, in addition to \eqref{level1ELrel}, we have \label{PropVWs2s4}	
		\begin{equation}
			V_3+W_3+2s_2t_4  = 0. \label{level1ELrel_2}
		\end{equation}
\end{Prop}					
\emph{Proof.} Restricting to $\mathcal{P}_s$, the generators $\theta_i$ of $\J_1$ satisfy congruences of the form:
				\[
					\ed\theta_i  \equiv -\pi_{i \alpha}\W \eta_\alpha +\tau_i, \mod \theta_1,\ldots,\theta_4 \quad (i = 1,\ldots,4).
				\]
			The tableau $(\pi_{i\alpha})$ now takes the form
				\[
					(\pi_{i\alpha}) = \left(
									      \begin{array}{cccccc}
									      		-\pi_1 & -\pi_2&0&0 &s_2\pi_8 &  s_2\pi_8\\
											-\pi_3&\pi_1-\pi_7 &0&0& -s_2\pi_7-\pi_9 & -\pi_9\\
											0&0&-\pi_4&-\pi_5 &-t_4\pi_{10} &t_4\pi_{10}\\
											0&0&-\pi_6 & \pi_7-\pi_4 & t_4\pi_7 +\pi_{11}& -\pi_{11}
									      \end{array}
									\right),
				\]	
			where, reduced modulo $\eta_0,\eta_1,\ldots,\eta_4$, 
				\begin{equation*}
					\begin{array}{llll}
					\pi_1 \equiv \varphi_1 - \psi_1, 
								&\pi_4 \equiv \varphi_5 - \psi_5,
								&	\pi_7 \equiv \varphi_0-\psi_0,
								&\pi_{10} \equiv \psi_6\\	[0.8em]
					\pi_2 \equiv \varphi_2 - \psi_2, 
								& \pi_5 \equiv \varphi_6+\psi_6,
								&	\pi_8 \equiv \psi_2,
								&\pi_{11} \equiv t_4\psi_5 -\ed(t_4).
								\\[0.8em]
					\pi_3 \equiv \varphi_3 - \psi_3, 
								&\pi_6 \equiv -\varphi_7-\psi_7,
								&\pi_9 \equiv  s_2 \psi_1-\ed(s_2),
								&			
					\end{array}
				\end{equation*}	
			By a calculation using Maple\texttrademark, it is easy to see that the torsion $\tau_i$ 
			can be absorbed only if the equations \eqref{level1ELrel} and \eqref{level1ELrel_2} hold.\qed
\vskip 2mm	
		
The equalities \eqref{level1ELrel} and \eqref{level1ELrel_2} tell us which Euler-Lagrange systems may be
special B\"acklund-related.
In particular, by Propositions \ref{lev1Perm} and \ref{VTrans}, it is easy to see that whether $V_2$ (hence $W_2$)
vanishes is independent of the choice of $1$-refined liftings.
Thus, we may locally\footnote{Namely, the conditions
									below hold on an entire open subset of $N$.}  
								classify special rank-$1$ B\"acklund transformations relating two hyperbolic Euler-Lagrange
systems into the following three types:
\begin{enumerate}[\qquad\qquad\qquad\bf Type I.]
	\item{$V_2 = 0$, $\det(S_1) = V_1V_4\ne 0$;}

	\item{$V_2 = 0$, $\det(S_1) = V_1V_4 = 0$;}

	\item{$V_2\ne 0$.}
\end{enumerate}	

\begin{Prop}
	$\bf(A)$ If a pair of hyperbolic Euler-Lagrange systems are related by a Type I special
	B\"acklund transformation, then one of them must be positive, the other negative.
\\
$\bf(B)$ Two hyperbolic Euler-Lagrange systems related by a Type III special B\"acklund transformation cannot be both degenerate. \label{SpecialBackLemma}
\end{Prop}

\emph{Proof.} For Part $\bf (A)$, it is immediate by \eqref{level1ELrel} that,
on a $1$-refined lifting,  $\det(S_1) = V_1V_4 = -W_1W_4 =  -\det(\bar S_1)\ne 0$.
Therefore, one of the Euler-Lagrange systems being B\"acklund-related is positive, the other negative.

To prove Part $\bf(B)$,
first apply Proposition \ref{lev1Perm} to show that, in 
this case, one can always find a $1$-refined lifting on which $V_1 = V_4 = 0$. 
(By Remark \ref{Hremark}, this can be achieved by acting on 
an initial $1$-refined lifting by an element $h\in H$ with $h_0 = h_1 = h_2 = 1$, $h_3 = V_4/V_2$, 
$h_4 = -V_1/V_2$.)
For such a $1$-refined lifting, by \eqref{level1ELrel_2}, the two Euler-Lagrange systems
can be both degenerate only when $s_2t_4 = 0$.  By Proposition \ref{s2t4=0}, both $S_1$ and $\bar S_1$
must vanish, which impossible since we have assumed $V_2 \ne 0$.
\qed
\vskip 3mm

Now we focus on Type $\rm II$.

Let $\phi: N\rightarrow M\times \bar M$ be a Type II special rank-$1$ B\"acklund transformation in the sense above. 
By Propositions \ref{VTrans} and \ref{lev1Perm}, 
whether the pair $(V_1,V_4)$ vanishes is 
independent of the choice of $1$-refined liftings of $\phi$. It follows that
$\phi$ must be one of the following two types.
			\begin{quote}
				\noindent	{\bf Type IIa}: on any $1$-refined lifting of $\phi$, $(V_1,V_4) = 0$;

				\noindent	{\bf Type IIb}: on any $1$-refined lifting of $\phi$, $(V_1,V_4) \ne 0$.		
			\end{quote}							
		\subsection{Type  IIa} 
		In this case, \eqref{level1ELrel} implies that 
	\[	
		W_1 = W_2 = W_4 = 0.
	\]	 	
	
				If locally either $V_3$ or $W_3$ is zero, which is independent of the choice of $1$-refined liftings, the 
				underlying B\"acklund transformation must relate a hyperbolic Euler-Lagrange system
				with the system corresponding to $z_{xy} = 0$. 
				See \cite{CI09} and \cite{zvyagin1991second} for a classification of all hyperbolic Monge-Amp\`ere systems that are 
				rank-$1$ B\"acklund-related to the equation $z_{xy} = 0$.
				
				More generally, we have the following theorem.
																
			\begin{theorem}\label{ThmIIa}
			If two hyperbolic Euler-Lagrange systems are related by a Type IIa special rank-$1$ 
			B\"acklund transformation, then each of them corresponds (up to contact equivalence) 
			to a second order PDE of the form 
								\[
									z_{xy} = F(x,y,z,z_x,z_y).
								\]						
			\end{theorem}									
\emph{Proof.} By definition, any Type IIa special B\"acklund transformation admits a $1$-refined lifting that is completely contained in the locus $\mathcal{P}_{{\rm IIa}}\subset \mathcal{P}_1$
						defined by the equations
								\[
									\mu = 1,\quad \epsilon = -1, \quad V_1 = V_2 = V_4 = W_1 = W_2 = W_4 = 0.
								\]
					  	Let $\mathcal{G}_2\subset \mathcal{G}_1$ be the subbundle defined by $V_1 = V_2 = V_4 = 0$; similarly, let $\bar{\mathcal{G}}_2\subset\bar{\mathcal{G}}_1$ be the subbundle defined 
						by $W_1 = W_2 = W_4 = 0$. It is clear that $\mathcal{P}_{{\rm IIa}}$ is the
						 product of $\mathcal{G}_2$, $\bar{\mathcal{G}}_2$, and a space of parameters with coordinates $(s_2, t_4)$.	
						 
						 The calculations below are performed using Maple\texttrademark.	
						
						By \eqref{VonFiber}, on $\mathcal{G}_2$, there exist functions $P_{ij}$ and $V_{3j}$ such that
								\begin{align*}
									\varphi_2& = P_{20}\alpha_0+P_{21}\alpha_1+\cdots+P_{24}\alpha_4,\\
									\varphi_6& = P_{60}\alpha_0+P_{61}\alpha_1+\cdots+P_{64}\alpha_4,\\
									\ed(V_3)& = (\varphi_1+\varphi_5)V_3 + V_{30}\alpha_0+V_{31}\alpha_1+\cdots+V_{34}\alpha_4.
								\end{align*}
						Similarly, on $\bar{\mathcal{G}}_2$, there exist functions $Q_{ij}$ and $W_{3j}$ such that
								\begin{align*}
									\psi_2& = Q_{20}\beta_0+W_{21}\beta_1+\cdots+Q_{24}\beta_4,\\
									\psi_6& = Q_{60}\beta_0+W_{61}\beta_1+\cdots+Q_{64}\beta_4,\\
									\ed(W_3)& = (\psi_1+\psi_5)W_3 + W_{30}\beta_0+W_{31}\beta_1+\cdots+W_{34}\beta_4.
								\end{align*}		
						There is freedom to add linear combinations of $\alpha_0,\ldots,\alpha_4$ (resp. $\beta_0,\ldots,\beta_4$) into $\varphi_i$ (resp. $\psi_i$) without changing the form of the corresponding 
						Monge-Amp\`ere structure equations. Using this, we can arrange the following 
						expressions to be zero:
								\[
									P_{21},~ P_{22},~ P_{63},~ P_{64};~ Q_{21},~ Q_{22},~ Q_{63},~ Q_{64}. 
								\]

						Applying $\ed^2 = 0$ to the Monge-Amp\`ere structure equations yields that
								\begin{equation*}
									\begin{alignedat}{2}
									\ed(\ed\alpha_1) &\equiv P_{24} V_3 \alpha_0\W \alpha_3\W \alpha_4&&\mod \alpha_1,\alpha_2,\\
									\ed(\ed \alpha_2) & \equiv V_{34} \alpha_0 \W \alpha_3 \W \alpha_4 &&\mod \alpha_1, \alpha_2,\\
									\ed(\ed \alpha_3) &\equiv P_{62} V_3 \alpha_0\W \alpha_1 \W \alpha_2 &&\mod \alpha_3, \alpha_4,\\
									\ed(\ed \alpha_4) &\equiv V_{32} \alpha_0\W \alpha_1\W \alpha_2 &&\mod \alpha_3, \alpha_4.
									\end{alignedat}
								\end{equation*}	
						This implies that, on $\mathcal{G}_2$,
								\begin{equation}
									P_{24}V_3 = V_{34} = P_{62}V_3 = V_{32} = 0.		\label{IIaP}
								\end{equation}
						By a similar argument, one can show that, on $\bar{\mathcal{G}}_2$,
								\[
									Q_{24}W_3 = W_{34} = Q_{62} W_3 = W_{32} = 0.	\label{IIaQ}
								\]		
						
						Restricted to $\mathcal{P}_{{\rm IIa}}$, the generators $\theta_i$ $(i = 1,\ldots,4)$ of $\J_1$ satisfy congruences of the form
								\[
									\ed\theta_i \equiv -\pi_{i\alpha}\W \eta_\alpha +\tau_i \mod\theta_1,\ldots,\theta_4.
								\]
						The tableau takes the form:	
								\[
									(\pi_{i\alpha}) = \left(
														\begin{array}{cccccc}
															-\pi_1 & 0&0&0&0&0\\
															-\pi_2 &\pi_1-\pi_5 &0&0&-s_2\pi_5+\pi_6 & \pi_6\\
															0&0&-\pi_3 &0&0&0\\
															0&0&\pi_4&-\pi_3+\pi_5 & t_4\pi_5-\pi_7 &\pi_7
														\end{array}
													\right),
								\]		
						where, modulo $\eta_i$, 
								\begin{equation*}
								\begin{array}{llll}
									\pi_1 \equiv \varphi_1 - \psi_1, &\pi_3 \equiv \varphi_5-\psi_5,
																& \pi_5 \equiv \varphi_0-\psi_0,
															&\pi_7 \equiv \ed(t_4)-t_4 \psi_5,\\[0.8em]
									\pi_2 \equiv \varphi_3 -\psi_3, &\pi_4 \equiv \varphi_7+\psi_7,
																&\pi_6 \equiv \ed(s_2) - s_2 \psi_1.
																&
								\end{array}								
								\end{equation*}
						Assuming $s_2, t_4$ to be both nonzero, 
						we compute and find that the torsion 
						can be absorbed only if the following expressions are zero:
								\[
									P_{20},~ P_{23}, ~P_{24}, ~P_{60}, ~P_{61}, ~P_{62};~ Q_{20}, ~Q_{23}, ~Q_{24}, ~Q_{60}, ~Q_{61}, ~Q_{62}.
								\]	
						One can verify that, on the subbundle of $\mathcal{G}_2$ defined by the vanishing of $P_{20}, P_{23}, P_{24},$ $ P_{60}, P_{61}$ and $P_{62}$, the following structure equations hold:
								\begin{align*}
									\ed\alpha_0 & = -\varphi_0\W \alpha_0 + \alpha_1\W \alpha_2 + \alpha_3\W \alpha_4,\\
									\ed\alpha_1& = -\varphi_1 \W \alpha_1,\\
									\ed\alpha_2& = -\varphi_3\W \alpha_1+(\varphi_1-\varphi_0)\W \alpha_2+V_3 \alpha_0\W \alpha_3,\\
									\ed\alpha_3& = -\varphi_5\W \alpha_3,\\
									\ed\alpha_4& = -\varphi_7\W \alpha_3 +(\varphi_5-\varphi_0)\W \alpha_4 +V_3 \alpha_0\W \alpha_1.
								\end{align*}					
						Clearly, the systems $\<\alpha_1\>$ and $\<\alpha_3\>$ are both integrable. 
						It is a similar case for the structure equations on $\bar{\mathcal{G}_2}$. 
						By Proposition \ref{zxyProp}, the proof is complete. \qed	
				
				\subsection{Type IIb}  
				
				In this case, on a $1$-refined lifting of $\phi$, either $V_1$ or $V_4$ vanishes. 
				By Proposition \ref{lev1Perm} (in particular, using $h_J$), we can
						 arrange $V_1\ne 0$ and $V_4 = 0$ on a $1$-refined lifting of $\phi$. 
				Such a $1$-refined lifting can always be chosen to further satisfy
				$V_1 = 1$ and $V_3 = W_3 = 0$ or $1$. 
				
				In the next proposition we show that 
						the case of $V_1 = 1$ and $V_3 = W_3 = 0$ is impossible. 
				Then we characterize the case when $\phi$ admits a $1$-refined lifting for which
						$V_1 = V_3 = W_3 = 1$.		
						
				\begin{Prop}Restricting to the locus in $\mathcal{P}_1$ defined by 
								\[	
									\mu = 1,~ \epsilon = -1,~V_1 = -W_1 = 1,~ V_2 = V_3 = V_4 = W_2 = W_3 = W_4 = 0, 
								\]
						    $\J_1$ has no integral manifold.		
				\end{Prop}
	\emph{Proof.} 
									  If there exists a $1$-refined lifting of a special B\"acklund transformation such that $V_3 = W_3 = 0$, then the equality \eqref{level1ELrel_2} enforces that 
									   $s_2t_4 = 0$ on such a lifting. By Proposition \ref{s2t4=0}, both Monge-Amp\`ere systems must be contact equivalent to the wave equation $z_{xy} = 0$.
									   In particular, $V_i$ and $W_i$ must all be zero on $\mathcal{G}_1$ and $\bar{\mathcal{G}_1}$, respectively. This contradicts our assumption. 
										\qed
					\vskip 3mm

					\begin{theorem} 
					Let $(M,\I)$ and $(\bar M,\bar\I)$ be two hyperbolic Euler-Lagrange systems.
					If $\phi:N\rightarrow M\times \bar M$ defines a Type IIb special rank-$1$ B\"acklund transformation relating $(M,\I)$ and $(\bar M,\bar \I)$, 
					then each of $(M,\I)$ and $(\bar M,\bar\I)$ must have 
					a characteristic system that contains a rank-1 integrable subsystem.	\label{PropIIb}
					\end{theorem}
					\emph{Proof.} The idea is similar to that of Theorem \ref{ThmIIa}. 
					We restrict the differential ideal $\J_1$ 
					to the locus $\mathcal{P}_{{\rm IIb}}\subset \mathcal{P}_1$ 
					defined by the equations
						\[
							\mu = 1,~ \epsilon = -1,~V_1 = -W_1 = V_3 = W_3 = 1,~ V_2 =  V_4 = W_2 =  W_4 = 0
						\]
					and analyze the obstructions to integrability of the resulting rank-$4$ Pfaffian system.
					
					The calculations below are performed using Maple\texttrademark.	
						
					By \eqref{VonFiber}, on the subbundle of $\mathcal{G}_1$
					defined by $V_1 = V_3  = 1$ and $V_2 = V_4 = 0$, there exist functions $P_{ij}$ such that 
						\begin{align*}
							\varphi_2 & = P_{20}\alpha_0+P_{21}\alpha_1+\cdots+P_{24}\alpha_4+\varphi_0  - \varphi_1+ \varphi_5,\\
							\varphi_3 & = P_{30}\alpha_0+P_{31}\alpha_1+\cdots+P_{34}\alpha_4+\varphi_1+\varphi_5,\\
							\varphi_6 & = P_{60}\alpha_0+P_{61}\alpha_1+\cdots+P_{64}\alpha_4.
						\end{align*}
					Using the freedom in the choice of $\varphi_i$, we can arrange that
					\[P_{21}, P_{22}, P_{23}, P_{31}, P_{32}, P_{34}, P_{64}\] are zero.

					Expanding $\ed(\ed\alpha_i) = 0$, we find that 
						\[
							P_{24} = 0, \quad P_{61} = -P_{62}, \quad P_{63} = 0.
						\] 	
					Similarly, on the subbundle of $\bar{\mathcal{G}}_1$
					defined by $-W_1 = W_3  = 1$ and $W_2 = W_4 = 0$, there exist functions $Q_{ij}$ such that 
						\begin{align*}
							\psi_2 & = Q_{20}\beta_0+Q_{21}\beta_1+\cdots+Q_{24}\beta_4-\psi_0 + \psi_1- \psi_5,\\
							\psi_3 & = Q_{30}\beta_0+Q_{31}\beta_1+\cdots+Q_{34}\beta_4-\psi_1-\psi_5,\\
							\psi_6 & = Q_{60}\beta_0+Q_{61}\beta_1+\cdots+Q_{64}\beta_4.
						\end{align*}
					Using the freedom in the choice of $\psi_i$, we can arrange that
					\[Q_{21}, Q_{22}, Q_{23}, Q_{31}, Q_{32}, Q_{34}, Q_{64}\] are zero. 
					
					By expanding $\ed(\ed\beta_i) = 0$, we find that
						\[
							Q_{24} = 0, \quad Q_{61} = Q_{62}, \quad Q_{63} = 0.
						\]
					
					Denote the restriction of ${\mathcal{J}}_1$ to  $\mathcal{P}_{{\rm IIb}}$ by ${\mathcal{J}}_{\rm IIb}$. By computation, we find that the torsion of $(\mathcal{P}_{{\rm IIb}},\J_{\rm IIb})$ 
					can be absorbed only if the following expressions are zero:
						\[
							s_2t_4+1, \quad P_{20} - Q_{20}, \quad P_{60}, \quad P_{62}, \quad Q_{60}, \quad Q_{62}.
						\]	
					In particular, the vanishing of $P_{60}, P_{62}, Q_{60}$ and $Q_{62}$ implies that
						\[
							\ed\alpha_3 = -\varphi_5\W \alpha_3, \qquad \ed\beta_3 = -\psi_5\W \beta_3.
						\]		
					The conclusion of the proposition follows. \qed

\section{New Examples of Type III}\label{newExamples}

In this section, we present some new examples of Type III special rank-$1$ B\"acklund transformations.
Their existence shows that a degenerate hyperbolic Euler-Lagrange system
may be special B\"acklund-related to a non-degenerate one. One of these examples (Section \ref{specialCase}) is, up to contact 
transformations,
a B\"acklund transformation relating solutions of the PDE
		\begin{equation}
			z_{xy} = \frac{2z}{(x+y)^2}	\label{GoursatPDE}
		\end{equation}
to those of a more complicated PDE of the form \eqref{MAequation} whose coefficients are given by
\eqref{ABCDE} in Appendix \ref{PDEForms}. We note that \eqref{GoursatPDE} 
is on the Goursat-Vessiot list\footnote{See \cite{goursat}, \cite{vessiot39} and \cite{vessiot42}.
The recent \cite{CI09} has a summary of the list.} of PDEs of the form
		\[
			z_{xy} = F(x,y,z,z_x,z_y)
		\]
that are Darboux-integrable at the 2-jet level.

\subsection{A Class of New Examples}

Using a method
in \cite{HuBacklund1} and after a somewhat lengthy calculation using Maple\texttrademark, we obtain on a 6-manifold $N$
involutive\footnote{Namely, exterior differentiation applied to these equations
yields identities.} structure equations (see Appendix \ref{appx1} for details)
of the form:
	\begin{equation}
		\ed\omega^i = -\frac{1}{2}C^i_{jk}(R,S,T)\omega^j\W\omega^k,	\label{primstrEqn}
	\end{equation}
	\begin{equation}\begin{split}
		\ed R &= R_{k}(R,S,T)\omega^k, \\ 
		\ed S &= S_{k}(R,S,T)\omega^k, \\
		\ed T &= T_{k}(R,S,T)\omega^k, 
	\end{split}	\label{augstrEqn}
	\end{equation}
where $C^i_{jk} = -C^i_{kj}, R_k, S_k, T_k$
are certain fixed analytic functions defined on $\Omega:=\{RST\ne 0\}\subset \R^3$.

One can verify that each of the ideals
	\[
		\I:=\<\omega^1,\omega^3\W\omega^4,\omega^5\W\omega^6\>_\alg, \quad \bar\I:=\<\omega^2,\omega^3\W\omega^4,\omega^5\W\omega^6\>_\alg
	\]
is generated by the pull-back of a hyperbolic Monge-Amp\`ere ideal ($\I_0$ and $\bar \I_0$, resp.) defined on
a 5-dimensional quotient ($M$ and $\bar M$, resp.) of $N$.

By
	\begin{equation*}
		\begin{alignedat}{2}
		\ed\omega^1&\equiv  - 
				R\omega^3\W\omega^4 + \omega^5\W\omega^6,&&\mod \omega^1,\omega^2,\\
		\ed\omega^2&\equiv \omega^3\W\omega^4 + \frac{1}{R}\omega^5\W\omega^6, &&	\mod \omega^1,\omega^2,
		\end{alignedat}
	\end{equation*}
it is immediate that $N$ is a special rank-$1$ B\"acklund transformation relating $(M,\I_0)$ and $(\bar M,\bar \I_0)$.
Furthermore, we choose the new bases
	\begin{equation}
		(\sigma^0,\sigma^1,\ldots,\sigma^4) = (\omega^1, -R\omega^3, \omega^4-\omega^1, \omega^5,\omega^6 + R^{-1}\omega^1) 	\label{sigmas}
	\end{equation}
and
	\begin{equation}
		(\tau^0,\tau^1,\ldots,\tau^4) =(\omega^2,\omega^3, \omega^4- R\omega^2,
						R^{-1}\omega^5,\omega^6-\omega^2),	\label{taus}
	\end{equation}
such that their pull-backs via arbitrary sections $\iota:M\rightarrow N$ and $\bar\iota:\bar M\rightarrow
N$, respectively, are $1$-adapted hyperbolic Monge-Amp\`ere coframings.

Computing using these new bases, we find that
	\[
		S_1 = \(\begin{array}{cc}
			2S& -R^{-1}T\\
			4T^{-1}S^2R&-2S
		\end{array}\), \quad S_2 = \bs0,
	\]
	\[	 \bar S_1 = \(\begin{array}{cc}
			-2RS & R^{-1}T\\
			2T^{-1}R(2S^2R^2 - T)& -2RS
		\end{array}\), \quad \bar{S_2} = \bs0.
	\]
Clearly, this tells us that both $\I_0$ and $\bar \I_0$ are hyperbolic Euler-Lagrange (Proposition \ref{EL}). Moreover,
since $\det(S_1)= 0$ and $\det(\bar{S_1}) = 2T\ne 0$, $\I_0$ is degenerate, and $\bar\I_0$ is non-degenerate,
by our definition.

By Cartan's third theorem (see \cite{Br14}), for each $u\in \R^6$ and for any
$(R_0,S_0,T_0)\in \Omega$, there exists a coframing  $(\omega^i)$ and functions 
$R,S,T$ on a neighborhood $U$ of 
$u$ such that the structure equations \eqref{primstrEqn} and \eqref{augstrEqn} hold
with $(R,S,T)(u) = (R_0,S_0,T_0)$.

This fact tells us that $T$ may be positive or negative. In other words, $\bar\I_0$ can be a hyperbolic
Euler-Lagrange system of either the positive or the negative type.

Note that, generically, the map $(R,S,T): N^6\rightarrow \R^3$ has rank $3$.
This is easy to see by computing the differentials $\ed(R), \ed(S)$ and $\ed(T)$. In this 
generic case, $N$
represents a B\"acklund transformation of cohomogeneity $3$. 

On the other hand, the map $(R,S,T)$ has rank $2$, which is the minimum rank possible, if and only if
		\[
			T = 2R^2S^2.
		\]
By the structure equations, 
$\ed(T - 2R^2S^2)$ vanishes whenever $T - 2R^2S^2$ does. 

This gives rise to the condition when   
the underlying B\"acklund transformation has cohomogeneity $2$, that is, when $T = 2R^2S^2$. 
In particular, since $T>0$ in this case,
$\bar \I_0$ must be of the positive type. We present further analysis of this case below.

\subsection{The $T = 2R^2S^2$ Case}\label{specialCase}\

All calculations below, unless otherwise noted, are performed using Maple\texttrademark.

Setting $T = 2R^2S^2$, the structure equations \eqref{primstrEqn}
and \eqref{augstrEqn} now take the form
\begin{equation}
		\ed\omega^i = -\frac{1}{2}\tilde C^i_{jk}(R,S)\omega^j\W\omega^k,	\label{CH2primstrEqn}
\end{equation}
{\begin{equation}\begin{split}
		\ed R &= \tilde R_{k}(R,S)\omega^k, \\ 
		\ed S &=\tilde S_{k}(R,S)\omega^k, \\
	\end{split}	\label{CH2augstrEqn}
	\end{equation}}
where $\tilde C^i_{jk}(R,S) = C^i_{jk}(R,S,2R^2S^2)$, and similarly for $\tilde R_k$ and $\tilde S_k$.

Let $\sigma^i$ and $\tau^i$ be as in \eqref{sigmas} and \eqref{taus}, respectively.
We find that
the $1$-form
	\[
		\phi:= -\frac{\omega^2}{2S} - \frac{1}{2SR}\sigma^0- \frac{1}{2}\sigma^1+SR\sigma^2
		-\frac{R}{2}\sigma^3 +\frac{R^2S^2 - S^2+1}{S}\sigma^4
	\]
is exact. Thus, locally there exists a function $f$ on $N$ such that
	\[
		\ed f = \phi.
	\]
Now consider the $1$-forms
	\[
		(\xi^0,\xi^1,\ldots,\xi^4) := \(\frac{\sigma^0}{e^f R},~ \frac{\sigma^1 - RS\sigma^2}{R},~
				\frac{\sigma^2}{e^f},~ \sigma^3 - RS\sigma^4, ~\frac{\sigma^4}{e^fR}\).
	\]
Using \eqref{CH2primstrEqn} and \eqref{CH2augstrEqn} and by letting $U = e^f/S$,  we obtain the following structure equations
for $(\xi^i)$:
{\small	
	\begin{align}
		\ed\xi^0 &= \xi^1\W\xi^2+\xi^3\W\xi^4,	\label{xi0}\\
	\ed\xi^1& = -U\xi^0\W\xi^1 + \frac{U}{2R}\xi^1\W\xi^2 + \frac{R^2-1}{2R}\xi^1\W\xi^3 
			-\frac{UR}{2}\xi^1\W\xi^4,		\\
	\ed\xi^2& = 2\xi^0\W\xi^1+U\xi^0\W\xi^2+2\xi^0\W\xi^3 - \frac{R^2+1}{2R}\xi^1\W\xi^2\\
			&+ R\xi^1\W\xi^4 - \frac{R^2-1}{2R}\xi^2\W\xi^3 +\frac{UR}{2}\xi^2\W\xi^4,	\nonumber\\
	\ed\xi^3& = U\xi^0\W\xi^3 - \frac{R^2-1}{2R}\xi^1\W\xi^3 + \frac{U}{2R}\xi^2\W\xi^3+
			\frac{UR}{2}\xi^3\W\xi^4,	\\
	\ed\xi^4&=2\xi^0\W\xi^1+2\xi^0\W\xi^3 - U\xi^0\W\xi^4 +\frac{R^2-1}{2R}\xi^1\W\xi^4\\
	&+\frac{1}{R}\xi^2\W\xi^3 -\frac{U}{2R}\xi^2\W\xi^4 + \frac{R^2+1}{2R}\xi^3\W\xi^4,	\nonumber\\
	\ed(U) &= -RU\xi^1+\frac{U}{R}\xi^3,\\
	\ed(R) &= UR\xi^0 +\frac{1}{2}(R^2+1)(\xi^1+\xi^3) +\frac{U}{2}\xi^2 - \frac{R^2U}{2}\xi^4.	\label{dR}
\end{align}
}

It is clear from these structure equations that $\xi^i$ $(i = 0,\ldots,4)$ are well-defined on a 
$5$-dimensional quotient $M^5$ of $N$, the corresponding hyperbolic Euler-Lagrange system 
being $(M,\I_0)$, where
		\[
			\I_0= \<\xi^0,\xi^1\W\xi^2,\xi^3\W\xi^4\>_\alg.
		\]
In particular, note that $\xi^1$ and $\xi^3$ are both integrable. By Proposition \ref{zxyProp},
$(M,\I_0)$ is contact equivalent to a PDE of the form $z_{xy} = F(x,y,z, z_x,z_y)$.

In Appendix \ref{PDEForms}, we find local coordinates for $(M,\I_0)$ and prove that
it is contact equivalent to the equation:
			\[
				z_{xy} =\frac{2z}{(x+y)^{2}} .
			\]

On the other hand, consider the $1$-forms
	\begin{align*}
		&(\eta^0, \eta^1,\ldots,\eta^4) :=\\ 
			&	\(RS\tau^0, ~R\tau^1+\frac{R(R^2S^2-1)}{S(R^2+1)}\tau^2,~
			S\tau^2, ~R\tau^3 - \frac{R(R^2S^2+1)}{S(R^2+1)}\tau^4, ~S\tau^4	
				\).
	\end{align*}
Using \eqref{CH2primstrEqn} and \eqref{CH2augstrEqn}, we find that
{\small	
\begin{align}
	\ed\eta^0 &= \frac{R^2+1}{2R^2}\eta^0\W\(\eta^1- \frac{R(2R^2-1)}{R^2+1}\eta^2 +R\eta^3 +\frac{R^2}{R^2+1}\eta^4\) 	\label{eta0}\\
			&+\eta^1\W\eta^2+\eta^3\W\eta^4,			\nonumber	\\
	\ed\eta^1& = -2\eta^0\W\(\frac{\eta^1}{R}+\eta^2+\eta^3+ \frac{R(R^2-1)}{R^2+1}\eta^4\)
	-\frac{R(2R^2+1)}{R^2+1}\eta^1\W\eta^2 \\
	&- \frac{1}{R}\eta^1\W\eta^3 + \frac{2}{R^2+1}\eta^1\W\eta^4 - \frac{R^2+2}{R^2+1}\eta^2\W\eta^3
	+\frac{2R}{(R^2+1)^2}\eta^2\W\eta^4,\nonumber
	\\
	\ed\eta^2& = \eta^0\W\(\frac{R^2+1}{R}\eta^2 - 2\eta^4\) - \frac{3R^2-1}{2R^2}\eta^1\W\eta^2
			 - \frac{1}{R}\eta^1\W\eta^4\\
			 &+\frac{R^2+3}{2R}\eta^2\W\eta^3 +\frac{R^2-3}{2(R^2+1)}\eta^2\W\eta^4,\nonumber\\
	\ed\eta^3 &= \eta^0\W\(2\eta^1 - \frac{2R(R^2-1)}{R^2+1}\eta^2 + \frac{R^2-1}{R}\eta^3\)
			+\frac{R^2-1}{2R^2}\eta^1\W\eta^3\\& - \frac{R}{R^2+1}\eta^1\W\eta^4
			+\frac{3R^2-1}{2R(R^2+1)}\eta^2\W\eta^3+\frac{R^2(R^2-1)}{(R^2+1)^2}\eta^2\W\eta^4
			\nonumber\\
			&+\frac{3R^2+1}{2(R^2+1)}\eta^3\W\eta^4,\nonumber\\
	\ed\eta^4&= 2\eta^0\W\eta^2 - \eta^1\W\eta^4 + \eta^2\W\eta^3+\frac{R(R^2-1)}{(R^2+1)}\eta^2\W\eta^4 - \frac{2}{R}\eta^3\W\eta^4,\\
	\ed(R) &= -(R^2+1)\eta^0 - \frac{R^2+1}{2R}\eta^1 - \frac{1}{2}\eta^2+\frac{R^2+1}{2}\eta^3
	-\frac{R}{2}\eta^4.	\label{anotherdR}
\end{align}
}

It is clear, by these structure equations, that $(\xi^i)$ descend to a coframing on a $5$-dimensional quotient
$\bar M^5$ of $N$, the corresponding hyperbolic Euler-Lagrange system 
being $(\bar M^5, \bar\I_0)$, where
	\[
		\bar\I_0 = \<\eta^0,\eta^1\W\eta^2,\eta^3\W\eta^4\>_\alg.
	\]
	
We have found local coordinates on $\bar M$ that put $(\bar M,\bar\I_0)$ (up to contact transformations) in a 
PDE form. Since both the expression and its derivation are rather complicated, we include them 
in Appendix \ref{PDEForms}.

{\remark This example shows the possibility for a pair of hyperbolic Euler-Lagrange systems of distinct types
to be special B\"acklund-related, 
in which one admits nontrivial first integrals for its characteristic systems, while the other
doesn't. 
In contrast, suppose that $S$ and $\bar S$ are, respectively, integral surfaces
of two hyperbolic Monge-Amp\`ere systems $(M,\I)$ and $(\bar M,\bar\I)$. If $S, \bar S$ are related to each other by a rank-$1$ B\"acklund transformation, then the characteristic curves in $S$ correspond, 
under the B\"acklund transformation, to the characteristic curves in $\bar S$, and \emph{vice versa}.
This is easy to see by the condition $(2)$ in Definition \ref{BacklundDef}.}

\section{Some Classical Examples}

In this section, we present some examples of special rank-$1$ B\"acklund transformations.
These examples are not new. We hope they can serve as a motivation for
the open questions described in the next section.
\begin{enumerate}[{\bf 1.}]

\item{\emph{A class of homogeneous rank-$1$ B\"acklund transformations.}
	In \cite{C01}, it is shown that, if, on $N^6$, there exists a local coframing 
	$(\theta_1,\theta_2, \omega^1,\ldots,\omega^4)$ satisfying
{	\begin{align*}
		\ed\theta_1 & = \theta_1\W(\omega^1+\omega^3) + \omega^1\W\omega^2+\omega^3\W\omega^4,\\
		\ed\theta_2& = -\theta_1\W(\omega^1+\omega^3)+ \omega^1\W\omega^2-\omega^3\W\omega^4,\\
		\ed\omega^1& = (B_1(\theta_1+\theta_2) -\sigma \omega^4)\W\omega^2+B_1\theta_1\W\theta_2,\\
		\ed\omega^2& = (-\sigma B_3(\theta_1+\theta_2)+ B_1^{-1}B_3\omega^4)\W\omega^1+\omega^3\W\omega^4,\\
		\ed\omega^3& = (-B_3(\theta_1 - \theta_2) - \sigma\omega^2)\W\omega^4 + B_3\theta_1\W\theta_2,\\
		\ed\omega^4& = (\sigma B_1(\theta_1-\theta_2)+B_3^{-1}B_1\omega^2)\W\omega^3+\omega^1\W\omega^2,
	\end{align*}
}	
where $B_i$ are constants and $\sigma = \pm 1$, then the systems
	\[
		\<\theta_1,\omega^1\W\omega^2,\omega^3\W\omega^4\>_{\alg}, 
		\quad \<\theta_2,\omega^1\W\omega^2,\omega^3\W\omega^4\>_{\alg}
	\]
are hyperbolic Monge-Amp\`ere systems well-defined on some $5$-dimensional quotients of $N$.
Here $N$ is a homogeneous rank-$1$ B\"acklund transformation relating the two Monge-Amp\`ere systems.

It is clear, by the equation (see Lemma \ref{epsilonmulemma})
					\[
						(\ed\theta_1 + \lambda \ed \theta_2)^2 \equiv 0 \mod \theta_1, \theta_2,
					\] 
that $\epsilon\mu^4 = -1$. Therefore, these B\"acklund transformations are special.	

In particular, when $\sigma = -1$, $B_1 = -(\tau+1)^2/2$ and $B_3 = (\tau^{-1}+1)^2/2$ for some constant
$\tau\ne 0$, one obtains the classical 
B\"acklund transformation between $K = -1$ surfaces in $\E^3$ when $\tau>0$ and
a B\"acklund transformation between $K = 1$ surfaces in $\E^{2,1}$ when $\tau<0$.
} 
\vskip 2mm
\item{Let $(x,y,u,p,q)$ and $(X,Y,v,P,Q)$ be coordinates on the spaces of two sine-Gordon systems:
	\begin{align*}
		\I &= \left\<\ed u - p\ed x - q\ed y, \left(\ed p - \frac{1}{2}\sin(2u)\ed y\right)\W\ed x\right\>,\\
		\bar\I & = \left \<\ed v - P\ed X - Q\ed Y, \(\ed P - \frac{1}{2}\sin(2v)\ed Y\)\W\ed X\right\>.
	\end{align*}
The locus $N$ in the product manifold defined by the equations
		\[
			X = x, \quad Y = y, \quad p = P+\lambda\sin(u+v),\quad Q  = -q+\lambda^{-1}\sin(u - v)
		\]
corresponds to the B\"acklund transformation \eqref{SGT}.	

The pull-backs to $N$ of $\theta =\ed u - p\ed x - q\ed y $ and $\bar\theta = \ed v - P\ed X - Q\ed Y$
	satisfy
		\[
			(\ed\theta+\ed\bar\theta)^2\equiv (\ed\theta - \ed\bar\theta)^2 \equiv 0\mod
			\theta,\bar\theta.
		\]
	It follows that the B\"acklund transformation \eqref{SGT} is special.	
}

\end{enumerate}

\section{Open Questions}\label{openquestions}

Several results concerning special B\"acklund transformations relating two hyperbolic
Euler-Lagrange systems can be summarized in the diagram (Figure \ref{diagram}) below.
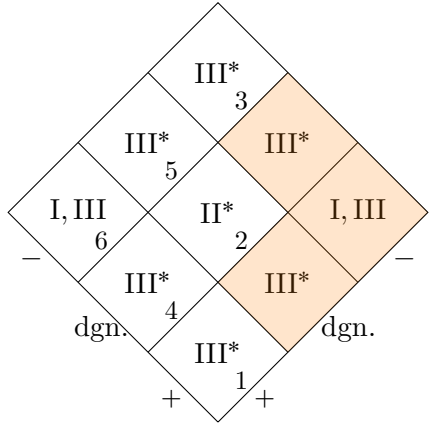
\begin{figure}[h!]
	\centering
		\begin{tikzpicture}[>=latex', scale=0.62]
			\draw[-] (0,0) --(-4.5, 4.5)--(0,9)--(4.5,4.5)--(0,0);
			\draw[-] (-1.5,1.5)--(3,6);
			\draw[-] (1.5,1.5)--(-3,6);
			\draw[-] (-3,3)--(1.5,7.5);
			\draw[-] (3,3)--(-1.5,7.5);
		
			\node at (-1,0.5) {$+$};
			\node at (-2.5,2) { dgn.~};
			\node at (-4, 3.5) {$-$};
			
			\node at (1,0.5) {$+$};
			\node at (2.7,2) {~dgn.};
			\node at (4, 3.5) {$-$};

			\node at (0, 7.5) {${\rm III}^*$};	
			\node at (0, 4.5) {${\rm II}^*$};		
			\node at (0, 1.5) {${\rm III}^*$};		
			
			\node at (-1.5, 3) {${\rm III}^*$};	
			\node at (1.5, 3) {${\rm III}^*$};	
			
			\node at (-3, 4.5) {${\rm I, III}$};	
			\node at (3, 4.5) {${\rm I, III}$};	
			
			\node at (-1.5, 6) {${\rm III}^*$};	
			\node at (1.5, 6) {${\rm III}^*$};	
				
			\draw [fill=orange, opacity=0.2]
       (1.5,1.5) -- (0,3) -- (1.5,4.5) -- (0, 6) --(1.5, 7.5) --(4.5,4.5)--cycle;
			
			\node at (0.5,0.9) {\small$1$};	
			\node at (0.5,3.9) {\small$2$};	
			\node at (0.5,6.9) {\small$3$};	
			
			\node at (-1,2.4) {\small$4$};	
			\node at (-1,5.4) {\small$5$};	
			
			\node at (-2.5,3.9) {\small$6$};	
				
		\end{tikzpicture}
		\caption{\footnotesize These are possible types of special rank-$1$
		B\"acklund transformations relating two hyperbolic Euler-Lagrange systems.
		The types of an Euler-Lagrange system (positive, degenerate, negative) are
		denoted by $+$, ${\rm dgn.}$, $-$, respectively.
		The superscript $^*$ indicates that example(s) exist and are either mentioned or 
		discovered in this paper.
		By symmetry, one only needs to consider regions $1$ to $6$ in the diagram.
		}\label{diagram}
\end{figure}
\begin{enumerate}[{\bf a.}]
	\item{An immediate question is: \emph{Can a special rank-$1$ B\"acklund transformation 
	 belong to region $6$ in the diagram?}}
	\item{Does there exist any special rank-$1$ B\"acklund transformation in regions
	$1$ and $3$ that does not relate an Euler-Lagrange system to itself?}
	\item{A more subtle problem is to understand the generality of existence in each of the regions. 
		For instance, it would be interesting to determine integers $k,\ell$ such that:
	\emph{The space of special rank-$1$ B\"acklund transformations relating two positive hyperbolic Euler-Lagrange systems
		is parametrized by $k$ functions of $\ell$ variables.}}
	\item{What can be said about the existence of non-special rank-$1$ B\"acklund transformations relating two hyperbolic Euler-Lagrange systems? Is there any hyperbolic Euler-Lagrange pair that are
	rank-$1$ B\"acklund-related but not special rank-$1$ B\"acklund-related?}
	\item{When can a hyperbolic Euler-Lagrange system be B\"acklund-related to
	a non-Euler-Lagrange hyperbolic Monge-Amp\`ere system?}
\end{enumerate}

\section{Acknowledgement}

	The author would like to thank his thesis advisor, Prof. Robert Bryant, for all his support
	and guidance in research and  Prof. Jeanne N. Clelland
	for her advice during the preparation of the current work for publication.
	He would  like to thank the referee for giving many 
	constructive comments, which motivated the writing of Appendix C.
	
\begin{appendix}
\section{Table of Structure Functions}\label{appx1}

In this appendix, we record the functions $C^i_{jk} = -C^i_{kj}, R_k, S_k, T_k$ involved in the
equations \eqref{primstrEqn} and \eqref{augstrEqn}.

\renewcommand*{\arraystretch}{1.8}
\begin{longtable}{||c|c||c|c||}
\caption{Expressions of $C^i_{jk}$, $R_k$, $S_k$, $T_k$ in \eqref{primstrEqn} and \eqref{augstrEqn}.
All $C^i_{jk}$ $(j<k)$ and $R_k, S_k, T_k$ that are not on this list are zero.}\\
\hline
$C^1_{12}$ & $ -\frac{SR^2}{T}$		& $C^1_{13}$ & $\frac{2R^2 - 1}{2R}$ \\
\hline
  	$C^1_{14}$&$\frac{R^2T+2R^2 - T}{4SR^3}$	&
$C^1_{15}$&$\frac{4R^2S^2+T}{2RT}$	\\
\hline
$C^1_{16}$&$-S$ 		&$C^1_{34}$ &$R$\\
\hline
$C^1_{56}$&$-1$ &	$C^2_{12}$		&	$\frac{4R^4S^2T+4R^4S^2 - R^2T^2-2R^2T+T^2}{4SR^3T}$	\\
\hline
$C^2_{23}$			&	$\frac{R(4R^2S^2 - 3T)}{2T}$ &
$C^2_{24}$ & $SR$ \\
\hline
 $C^2_{25}$ & $-\frac{R^2+2}{2R}$ & $C^2_{26}$ &$ \frac{R^2T+2R^2-T}{4SR^2}$\\
\hline
$C^2_{34}$ & $-1$ & $C^{2}_{56}$ &$-\frac{1}{R}$\\
\hline
 $C^3_{12}$&$-1$&	$C^{3}_{13}$& $\frac{2R^4S^2T+4R^4S^2- 2R^2S^2T-R^2T^2-2R^2T+T^2}{2SR^3T}$\\
 \hline
  $C^{3}_{14}$& $\frac{1}{R}$ & $C^3_{24}$ & $-1$\\
\hline
$C^3_{35}$ & $\frac{4R^2S^2 - R^2T - T}{2RT}$	& $C^3_{34}$ & $\frac{4R^4S^2 - R^2T-2R^2+T}{4SR^3}$
		\\
\hline
$C^3_{45}$ & $\frac{2S}{R}$&
$C^{3}_{46}$ & $-\frac{T}{R^2}$ \\
\hline
 $C^3_{36}$ & $-\frac{4R^2S^2 - R^2T - 2R^2+T}{4SR^2}$ & $C^4_{13}$&
$\frac{2(2R^2S^2 - T)}{RT}$\\
\hline
$C^4_{24}$ & $\frac{SR^2}{T}$ & $C^4_{14}$ & $-\frac{4R^4S^2-4R^2S^2T - R^2T^2-2R^2T+T^2}{4SR^3T}$
\\
\hline
 $C^4_{34}$ & $-\frac{1}{2R}$ &
$C^4_{36}$ & $1$ \\
\hline
 $C^{4}_{45}$ & $-\frac{4R^2S^2 -T}{2RT}$ & $C^{4}_{46}$ & $S$\\
\hline
$C^4_{56}$ & $-1$ &$C^5_{15}$&$-\frac{2R^4S^2T+4R^4S^2-2R^2S^2T-R^2T^2-2R^2T+T^2}{2SR^3T}$\\
\hline
$C^5_{12}$& $-1$&
$C^5_{16}$ & $1$ \\
\hline
 $C^5_{26}$ & $R$ & $C^5_{36}$ & $2SR^2$\\
\hline
$C^5_{46}$ & $T$ & $C^5_{45}$ & $-\frac{4R^4S^2+R^2T+2R^2-T}{4SR^3}$ \\
\hline
 $C^5_{35}$ & $-\frac{4R^4S^2-R^2T-T}{2RT}$ &
$C^5_{56}$ & $-\frac{4R^2S^2+R^2T+2R^2-T}{4SR^2}$ \\
\hline
$C^6_{15}$& $-\frac{2(2R^2S^2 - T)}{T}$&$C^6_{16}$ & $\frac{4R^4S^2T+4R^4S^2-R^2T^2-2R^2T+T^2}{4SR^3T}$\\
\hline
$C^6_{26}$ & $\frac{SR^2}{T}$ &$C^6_{34}$ & $-1$ \\
\hline
 $C^6_{36}$ & $\frac{R(4R^2S^2-T)}{2T}$&
$C^6_{45}$ & $-1$\\
\hline
$C^6_{46}$ & $SR$ & $C^6_{56}$ & $-\frac{R}{2}$\\
\hline
$R_3$ &$-\frac{4R^4S^2-R^2T+T}{2T}$ & $R_1$ &$\frac{2R^4S^2T+4R^4S^2-2R^2S^2T - R^2T^2-2R^2T+T^2}{2SR^2T}$\\
\hline
$R_4$ & $- \frac{4R^4S^2-R^2T-2R^2+T}{4SR^2}$ & $R_5$ & $\frac{4R^2S^2+R^2T-T}{2T}$\\
\hline
$R_6$ &$ - \frac{4R^2S^2+R^2T+2R^2-T}{4SR}$ &$S_1$&  $\frac{4R^4S^2T-4R^4S^2+8R^2S^2T-R^2T^2+2R^2T-3T^2}{4R^3T}$\\
\hline
$S_2$ & $-\frac{S^2R^2}{T}$ & $S_3$ &$ - \frac{SR(4R^2S^2- T)}{2T}$\\
\hline
$S_4$ & $- \frac{2R^2S^2-T}{2R}$ & $S_5$ & $-\frac{S(8R^2S^2+R^2T-2T)}{2RT}$\\
\hline
$S_6$ & $\frac{8R^2S^2+R^2T+2R^2-3T}{4R^2}$ &$T_1$&$\frac{4R^4S^2T+4R^4S^2-R^2T^2-2R^2T+T^2}{2SR^3}$\\
\hline
$T_2$ & $-2SR^2$ & $T_3$ & $- \frac{4R^4S^2+T}{R}$\\
\hline
$T_4$ & $-\frac{T(4R^4S^2-R^2T-2R^2+T)}{2SR^3}$ & $T_5$ & $-\frac{T}{R}$\\
\hline
\multicolumn{4}{r}{\textit{(End of Table 1)}} \\
\end{longtable}

\section{ $(M, \mathcal{I}_0)$ and $(\bar M,\bar{\mathcal{I}}_0)$ in PDE Forms} \label{PDEForms}

In this Appendix, we integrate the structure equations \eqref{xi0}-\eqref{dR} for $(M,\I_0)$
and \eqref{eta0}-\eqref{anotherdR} for $(\bar M,\bar\I_0)$. 
By `integrate', we mean finding local coordinates, expressing the coframings and structure functions
in terms of these coordinates such that the structure equations hold identically.
As a result, we put
these hyperbolic Euler-Lagrange systems in PDE forms (up to contact transformations).
Most calculations below are performed using Maple\texttrademark.
\vskip 3mm

\noindent{\bf I. Integration of $(M,\I_0)$.}

Consider \eqref{xi0}-\eqref{dR}. We integrate these structure equations in the following steps.


{
\begin{enumerate}[1)]

\item{First note that these structure equations are invariant under the transformation:
	\[
		(\xi^0,\xi^1,\xi^3)\mapsto(-\xi^0,-\xi^1,-\xi^3), \qquad (U,R)\mapsto (-U,-R).
	\]
As a result, we can assume that $R>0$.}

\item{It is clear that
	\[
		\ed\xi^i \equiv 0 \mod \xi^i \qquad (i = 1,3).
	\]
Thus, both $\xi^1$ and $\xi^3$ are multiples of exact forms.
In fact, we find that the $1$-forms
	\begin{align*}
		\psi_1	& = -U\xi^0 - \frac{3R^2+1}{2R} \xi^1 - \frac{U}{2R}\xi^2 -\frac{R^2-1}{2R}\xi^3+\frac{UR}{2}\xi^4,\\
		\psi_3 & = U\xi^0 - \frac{R^2 - 1}{2R}\xi^1 + \frac{U}{2R}\xi^2 +\frac{R^2+3}{2R}\xi^3 - \frac{UR}{2}\xi^4
	\end{align*}
satisfy
	\[
		\ed \xi^1 = \psi_1\W\xi^1, \quad \ed\xi^3 = \psi_3\W\xi^3
	\]
and
	\[
		\ed \psi_1 = \ed \psi_3 = 0.
	\]
It follows that there exist local coordinates $x,y,r,s$ such that
	\[
		\xi^1 = e^r\ed x, \quad \xi^3 = e^s\ed y; \quad \ed r = \psi_1, \quad \ed s = \psi_3,
	\]	
with freedom to add constants to $r,s,x,y$, simultaneously scaling $x,y$.
}

\item{It is straightforward to verify that 
	\[
		\ed(e^{r+s}/U^2) = \ed(e^rR/U) = 0.
 	\]	
Using the freedom of adding a constant to $r$ or $s$, we can arrange that
	\[
		e^{r+s} = U^2;
	\]
moreover, there exists a constant $\lambda\ne 0$ such that $e^rR = \lambda U$. 
(At this point, we choose not to normalize $\lambda$.)
Using these and the expression of $\ed U$, we obtain 
	\[
		\ed U = -U^2\lambda  \ed x+\frac{U^2}{\lambda} \ed y.
	\]
Direct integration gives
	\[
		U = \frac{1}{\lambda x - \lambda^{-1}y +C}
	\]	
for some constant $C$. We can normalize $C$ to be zero by adding a constant to $x$ or $y$.	
}

\item{Since
		\[
			\ed\xi^i\equiv 0\mod \xi^0,\xi^2,\xi^4\qquad (i = 0,2,4),
		\]
the system $\<\xi^0,\xi^2,\xi^4\>$ is Frobenius. We proceed by looking for its first integrals.	

Indeed, we find that 
	\[
		\Phi:= \xi^0 +\frac{1}{2R}\xi^2 -\frac{R}{2}\xi^4
	\]
satisfies
	\[
		\ed\Phi = \frac{R^2-1}{R}(\xi^1+\xi^3)\W\Phi.
	\]
Hence, $\Phi$ is a multiple of an exact form. By adding an appropriate multiple of $\Phi$ into $\frac{R^2-1}{R}(\xi^1+\xi^3)$,
we find that
	\[
		\ed(-2\ln R + 2\ln(R^2+1)) \equiv \frac{R^2-1}{R}(\xi^1+\xi^3)\mod \Phi.
	\]
Therefore, locally there exists a function $v$ such that 	
	\[
		\Phi = \(\frac{R^2+1}{R}\)^2 \ed v.
	\]
}
\item{ One can verify that the functions $R,x,y,v$ do not have linearly independent differentials.
	In fact, we find that
		\[
			\ed\(\frac{1}{R^2+1}\) = -\frac{\lambda U}{R^2+1}\ed x - \frac{R^2 U}{R^2+1}\ed y
					 - 2U \ed v.
		\]
	Setting $G(x,y,z):= (R^2+1)^{-1}$, direct integration gives
		\[
			G = \frac{1}{R^2+1} = -2Uv + \frac{-y+C_1}{\lambda^2x - y},
		\]	
	where $C_1$ is a constant. We can normalize $C_1$ to be zero by adding constants to $x$ and $y$.
	For convenience, we express $v$ in terms of $R,x,y$.
}
\item{Next, we observe that the $1$-forms
	\[
		\Xi_2:= \frac{1}{R}\xi^2 +\frac{\lambda}{R}\ed x - \frac{2\ln R}{\lambda}\ed y,
		\qquad
		\Xi_4 := R\xi^4 - 2\lambda \ln R~\ed x - \frac{R^2}{\lambda}\ed y
	\]
satisfy
	\[
		\ed(\Xi_2+\Xi_4) = 0,\qquad
		 \ed((\lambda^2 x -y)\Xi_2) = -\ed y \W (\Xi_2+\Xi_4).
	\]
Consequently, there exist functions $g,h$ such that
	\[
		 \Xi_2+\Xi_4 = \ed g,   \qquad \Xi_2 = \frac{\ed h - y\ed g}{\lambda^2 x - y}. 
	\]	
}

\item{	
At this point, each $\xi^i$ as well as the function $U$ can be expressed in terms of $x,y,R,g,h$ and their exterior derivatives.		
Moreover, \eqref{xi0}-\eqref{dR} become identities.

In particular, one can already put $\xi^0$ in the form
	\[
		\xi^0 = \ed z - p \ed x - q\ed y,
	\]
by introducing certain functions $z,p,q$. In our choice,
\[
z = \frac{\lambda^2x - y}{\lambda}\ln R +\frac{(\lambda^2x+y)g}{2(\lambda^2x - y)} - \frac{h}{\lambda^2x- y}.
\]

Computing using these coordinates, we obtain
	\[
		\xi^1\W\xi^2 + \ed p\W\ed x = \frac{2\lambda^2z}{(\lambda^2x - y)^2}\ed x\W\ed y.
	\]
This tells us that $(M,\I_0)$ is contact equivalent to the PDE
	\[
		z_{xy} = - \frac{2\lambda^2z}{(\lambda^2x - y)^2},
	\]	
which is, after we scale $x$ and $y$, contact equivalent to
	\[
		z_{xy} = \frac{2z}{(x+y)^2}.
	\]	
	
}

\end{enumerate}
}

\noindent{\bf II. Integration of $(\bar M,\bar\I_0)$.}

We integrate the structure equations \eqref{eta0}-\eqref{anotherdR} in the following steps.

{
\begin{enumerate}[1)]

\item{First note that 
	\[
		\ed\eta^i \equiv 0  \mod\eta^2,\eta^4\qquad (i = 2,4).
	\]
Thus, the system $\<\eta^2,\eta^4\>$ is Frobenius. In particular, there exists a linear combination
of $\eta^2,\eta^4$ that is a multiple of an exact $1$-form. Setting undetermined 
weights that are functions of $R$, we find that the $1$-form
	\[
		\Phi := \eta^2- R\eta^4
	\]
satisfies
	\[
		\ed\Phi = \psi\W\Phi
	\]	
for some $1$-form $\psi$.

By adding an appropriate multiple of $\Phi$ to $\psi$, we obtain an exact $1$-form (still denoted by $\psi$):
	\[
	\begin{split}
		\psi &= -\frac{R^2 - 1}{R}\eta^0 - \frac{3R^2-1}{2R^2}\eta^1 + \frac{R^2-3}{2R}\eta^3
			+\frac{(R^2-1)^2}{R^2+1}\eta^4 \\
			&+ \frac{2R^4 - R^2 + 1}{2R(R^2+1)}\Phi.
	\end{split}		
	\]	
Thus, locally there exist functions $f, g$ such that
	\[
		\ed f = \psi, \qquad \Phi = e^f \ed g.
	\]	
}
\item{Next, we compute and observe that
	\[
		\ed \eta^4\equiv \(-\frac{2(2R^2-1)}{R(R^2+1)}\ed R+ 2\ed f\)\W\eta^4\mod \ed g.
	\]
This motivates the calculation:
	\[
		\ed\(\frac{(R^2+1)^3}{e^{2f}R^2}\eta^4\) = \gamma\W\ed g,
	\]	
where
\[
	\begin{split}
		\gamma & = \frac{(R^2+1)^3}{R^4e^f}\eta^3+\frac{(R^2+1)^2}{Re^f}\eta^4
					 - \frac{(R^2+1)^2(3R^2-1)}{R^4e^f}\ed R\\
			& +\frac{(R^2+1)^3}{R^3e^f}\ed f+(R^2+4\ln|R| -R^{-2})\ed g.
	\end{split}
\]	
We constructed this $\gamma$, by choosing an appropriate `$\ed g$'-term, so that it is exact.

Consequently, there exist functions $u, v$ such that
	\[
		\ed u = \gamma, \qquad \ed v = \frac{(R^2+1)^3}{e^{2f}R^2}\eta^4 - u\ed g.
	\]
}

\item{Now each $\eta^i$ can be expressed completely in terms of $f, g, R, u, v$ and their differentials.
Moreover, \eqref{eta0}-\eqref{anotherdR} become identities.

In particular,
\[
	\begin{split}
	\eta^0& = \frac{3R^2-1}{R^2+1}\ed R - \frac{R^2e^f((R^2+1)^2+2(R^2+1)\ln|R|+ e^f uR)}{(R^2+1)^3}\ed g\\
	&+\frac{R^2e^f}{2(R^2+1)^2}\ed u - \frac{e^{2f}R^3}{(R^2+1)^3}\ed v -  \frac{R}{2}\ed f.
	\end{split}
\]
By the Pfaff theorem, up to scaling, we can put $\eta^0$ in the form 
\[
\ed z - p\ed x - q\ed y.
\]
We approach this by following the proof of the Pfaff theorem in Chapter II of \cite{BCG}.
We first expand $\ed \eta^0\W\eta^0\W\ed R$ in coordinates; from its expression we find 
that
\[
	x := 2(R^2+2\ln|R|+1)g - u
\]
satisfies
\[
	\ed \eta^0 \W\eta^0\W\ed R\W \ed x = 0.
\]
Next, writing $u$ in terms of $x, R$ and $g$, the $3$-form $\eta^0\W\ed R\W\ed x$ has only $3$
terms in it. This allows us to find 
\[
	y:= -v - (R^2+1+2\ln|R|)g^2  + xg +\frac{1}{4R^2e^{2f}}(R^2+1)^3,
\]
which satisfies
\[
	\eta^0\W\ed R\W\ed x\W\ed y = 0.
\]
}
\item{It follows from the previous step that $\eta^0$ is a linear combination of $\ed x, \ed y$ and $\ed R$.
After scaling $\eta^0$ such that the coefficient of $\ed R$ is $1$, we obtain the $1$-form
	\[
		\ed R - p \ed x - q\ed y,
	\]
where
\[
	p = \frac{(2gRe^f +R^2+1)Re^f}{(R^2+1)G}, \qquad q = -\frac{2R^2e^{2f}}{(R^2+1)G},
\]
and
\[
	G = 4ge^{f} (gRe^f + R^2 + 1) + R^3+R.
\]
}
\item{It is now reasonable to set $z = R$.
 	Using these coordinates, we compute and observe that
	\[\begin{split}
		\eta^1\W\eta^2&\equiv A\ed p\W\ed q+B\ed p\W\ed y+ C(\ed x\W\ed p - \ed y\W \ed q)\\
			&+ D\ed x\W\ed q + E \ed x \W \ed y \qquad \mod \eta^0,\ed \eta^0,
	\end{split}
	\]
where
\begin{equation}\begin{split}
	A& = 2qz(z^2+1)^3,\\
	B & = 2q^2(z^2+1)^2(4p^2z^3-qz^2+4p^2z+3q),\\
	C& = -2pq(z^2+1)^3(4p^2z+q),\\
	D& = (z^2+1)(4p^2z^3+ qz^2+4p^2 z - q)(2p^2z^2+qz+2p^2),\\
	E& = -q^3(4p^2z^5 +qz^4- 16p^2z^3-8qz^2-20p^2z-q).
	\end{split}	\label{ABCDE}
\end{equation}	

The PDE form of $(\bar M,\bar\I_0)$ (up to contact transformations) 
is therefore \eqref{MAequation} in the introduction, 
where $A,B,C,D,E$ are the same as the above with
$p,q$ being replaced by $z_x$ and $z_y$, respectively. 

Furthermore, hyperbolicity demands that 
\[	
	0<AE-  BD+C^2 = 8q^4(2p^2z^2+2p^2+zq)(z^2+1)^4;
\]
in other words, the domain of the variables $(x,y,z,p,q)$ needs to satisfy
 $q\ne 0$ and $2p^2z^2+2p^2+zq>0$.
}

\end{enumerate}
}

\section{A Note on Calculation} \label{CalcNote}

This appendix is written for those readers who wish to see certain concepts and examples in coordinates. In item {\bf I}, we provide a program that computes the hyperbolic
Monge-Amp\`ere relative invariants $S_1$ and $S_2$ in coordinates.
In items {\bf II}, {\bf III}, we focus on hyperbolic Euler-Lagrange systems and their types.
\\

\noindent{\bf I. The Monge-Amp\`ere Relative Invariants \eqref{S1S2def}.}

The Monge-Amp\`ere equation \eqref{MAequation} corresponds to the exterior differential 
system $(J^1(\R^2,\R),  \<\theta,\ed\theta, \Omega\>_\alg)$ with
	\begin{align*}
		\theta &= \ed z - p \ed x - q \ed y, \\
		 \Omega &= A\ed p \W \ed q+ B\ed p \W\ed y + C(\ed x\W \ed p - \ed y\W\ed q) \\
		 &+ D \ed x \W \ed q + E\ed x \W \ed y.
	\end{align*}
	
If $A,B,D,E$ are all zero, then the system is either empty or
 equivalent to the wave equation $z_{xy} = 0$.
From now on, we assume that not all of $A,B,D,E$ are zero. Under this assumption,
we could always make one of the following contact transformations to arrange that $E\ne 0$:
\begin{enumerate}[\quad 1)]
	\item{$(x',y',z',p',q')=(p,q,z-px - qy,-x,-y)$;}
	\item{$(x',y',z',p',q')=(p,y,z - px, -x, q)$;}
	\item{$(x',y',z',p',q')=(x,q,z - qy, x, -y)$.}
\end{enumerate}
It then follows that, by scaling the equation, we can arrange that $E = 1$.
In doing this, the hyperbolicity assumption would allow us to express 
\[
	A = BD - C^2 + \mu^2
\]
for some $\mu>0$.

By introducing the $1$-forms:
	\begin{align*}
		\eta^1 &= (C+\mu)\ed p + D \ed q + \ed y,\\
		\eta^2 &= -B \ed p +(\mu - C)\ed q - \ed x,\\
		\eta^3& = (C - \mu)\ed p + D\ed q + \ed y,\\
		\eta^4& = B\ed p + (\mu+C)\ed q + \ed x,
	\end{align*}
we obtain:
	\begin{equation*}
		\Omega+ \mu\ed \theta = \eta^1\W\eta^2,\qquad
		-\Omega + \mu\ed \theta = \eta^3\W\eta^4.
	\end{equation*}
In other words, 
	\[
		\<\theta, \ed\theta,\Omega\>_\alg = \<\theta,\eta^1\W\eta^2,\eta^3\W\eta^4\>_\alg.
	\]
Letting $\eta^0 = 2\mu\theta$, we have
	\[
		\ed\eta^0\equiv \eta^1\W\eta^2 +\eta^3\W\eta^4\mod\eta^0.
	\]
However, the coframing $(\eta^0,\eta^1,\ldots,\eta^4)$ is not 
necessarily $1$-adapted (see Section \ref{InvMA}).

Now define the coefficients $T^i_{jk} = -T^i_{kj}$
by
\[
	\ed\eta^i = \frac{1}{2}T^i_{jk}\eta^j\W\eta^k.
\]
In particular, let\footnote{We remark that each $c^i$ has 21 terms in $B,C,D,\mu$,
their partial derivatives, and $p,q$.} 
\[
	c^1 = T^1_{34},\quad c^2 = T^2_{34}, \quad c^3 = T^3_{12}, \quad c^4 = T^4_{12}.
\]
Construct a new coframing $(\omega^0,\omega^1,\ldots,\omega^4)$ by
\[
	\omega^0 = \eta^0, \quad 
	\omega^i = \eta^i - c^i\eta^0,\quad  (i = 1,\ldots,4). 
\]
One can verify that $(\omega^0,\omega^1,\ldots,\omega^4)$ is $1$-adapted; thus,
it can be used to compute the Monge-Amp\`ere relative invariants $V_1,\ldots,V_8$.

In fact, defining $\tilde T^i_{jk} = -\tilde T^i_{kj}$ using
\[
	\ed\omega^i = \frac{1}{2}\tilde T^i_{jk}\omega^j\W\omega^k,
\]
we obtain
\[
	S_1 =\(\begin{array}{cc}
			V_1	& V_2\\
			V_3 & V_4
			\end{array}\) =  \frac{1}{2}\(\begin{array}{cc}
			\tT^1_{03} - \tT^4_{02}	& \tT^1_{04} + \tT^3_{02}\\
			\tT^2_{03} + \tT^4_{01} & \tT^2_{04} - \tT^3_{01}
			\end{array}\), 
\]			
\[
	S_2 =\(\begin{array}{cc}
			V_5	& V_6\\
			V_7 & V_8
			\end{array}\)   = \frac{1}{2}\(\begin{array}{cc}
			\tT^1_{03} + \tT^4_{02}	& \tT^1_{04} - \tT^3_{02}\\
			\tT^2_{03} - \tT^4_{01} & \tT^2_{04} + \tT^3_{01}
			\end{array}\).	
\]
Such expressions of $V_i$ depend on up to the second partial
derivatives of $B,C,D,\mu$ and are rather complicated; to give the reader a sense, in our 
calculation, $V_1,\ldots,V_8$
 have 476, 159, 159, 476, 155, 262, 155 and 262 terms in them, respectively. 
We have used the following code in our calculation
(with Maple\texttrademark and the {\fontfamily{pcr}\selectfont Cartan}
package), where {\fontfamily{pcr}\selectfont M} stands for the function $\mu$ above.
\vskip 3mm

{\fontfamily{pcr}\selectfont \scriptsize
\setlength\parindent{0pt}
> restart: 

> with(Cartan):

> unprotect(D):

> Omega:= A*d(p)\mw d(q) + B*d(p)\mw d(y) + C*(d(x)\mw d(p)- d(y)\mw d(q)) 

\hspace{\codesp}+ D*d(x)\mw d(q) + E*d(x)\mw d(y):

> Theta:= d(x)\mw d(p) + d(y)\mw d(q):

> E:=1:

> A:= B*D - C$\W$2 + M$\W$2:

> var[1]:= x: var[2]:= y: var[3]:= z: var[4]:= p: var[5]:= q:

> for t in [B,C,D,M] do:

\hspace{\codesp}    d(t):= add(t[i]*d(var[i]), i = 1..5):
    
\hspace{\codesp}od:

> for t in [B,C,D,M] do:

\hspace{\codesp}    for i from 1 to 5 do:
    
\hspace{\codesp}        d(t[i]):= add(t[i,j]*d(var[j]), j = 1..5):
        
\hspace{\codesp}    od:
    
\hspace{\codesp}od:

> for t in [B,C,D,M] do:

\hspace{\codesp}   d(d(t)):Simf(\%):ScalarForm(\%):solve(\%):assign(\%):
   
\hspace{\codesp}od:

> for i from 0 to 4 do:

\hspace{\codesp}  Form(eta[i] = 1):

\hspace{\codesp}od:

> eta2coord:= \{eta[0] = (d(z)- p*d(x) - q*d(y))*(2*M),

\hspace{\codesp}eta[1] = (M + C)*d(p) + D*d(q)  + d(y),

\hspace{\codesp}eta[2] = -B*d(p) + (M - C)*d(q)  - d(x),

\hspace{\codesp}eta[3] = (-M + C)*d(p) + D*d(q)  + d(y), 

\hspace{\codesp}eta[4] = B*d(p) + (M + C)*d(q)  + d(x)\}:

> coord2eta:= solve(eta2coord, \{d(z), d(x), d(y), d(p), d(q)\}):

> for i from 0 to 4 do:

\hspace{\codesp}         d(subs(eta2coord, eta[i])):Simf(\%):subs(coord2eta,\%):deta[i]:=Simf(\%):

\hspace{\codesp}   for j from 0 to 3 do:
   
\hspace{\codesp}     for k from j+1 to 4 do:
         
\hspace{\codesp}         T[i,j,k]:= pick(deta[i], eta[j], eta[k]):
         
\hspace{\codesp}     od:
     
\hspace{\codesp}   od:
   
\hspace{\codesp}od:

> c1:= T[1,3,4]:

\hspace{\codesp}c2:= T[2,3,4]:

\hspace{\codesp}c3:= T[3,1,2]:

\hspace{\codesp}c4:= T[4,1,2]:

> for i from 0 to 4 do:

\hspace{\codesp}  Form(omega[i] = 1):

\hspace{\codesp}od:

> assign(eta2coord):

> om2coord:= \{omega[0] = eta[0],

\hspace{\codesp}omega[1] = eta[1] - c1*eta[0],

\hspace{\codesp}omega[2] = eta[2] - c2*eta[0],

\hspace{\codesp}omega[3] = eta[3] - c3*eta[0], 

\hspace{\codesp}omega[4] = eta[4] - c4*eta[0]\}:

> coord2om:= solve(om2coord, \{d(z), d(x), d(y), d(p), d(q)\}):

> for i from 0 to 4 do:

 \hspace{5mm}        d(subs(om2coord, omega[i])):Simf(\%):subs(coord2om,\%): dom[i]:= Simf(\%):
 
\hspace{\codesp}   for j from 0 to 3 do:
   
 \hspace{\codesp}    for k from j+1 to 4 do:
     
   \hspace{\codesp}      T[i,j,k]:= pick(dom[i], omega[j], omega[k]):
         
 \hspace{\codesp}    od:
     
\hspace{\codesp}   od:
   
\hspace{\codesp}od:

> V[1]:= (T[1,0,3] - T[4,0,2])/2:
 	V[2]:= (T[1,0,4] +  T[3,0,2])/2:
 
\hspace{\codesp}V[3]:= (T[2,0,3] + T[4,0,1])/2:
			  V[4]:= (T[2,0,4] - T[3,0,1])/2:
  
\hspace{\codesp}V[5]:= (T[1,0,3] + T[4,0,2])/2:
 	V[6]:= (T[1,0,4] - T[3,0,2])/2:
 
\hspace{\codesp}V[7]:= (T[2,0,3] - T[4,0,1])/2:
	 V[8]:= (T[2,0,4] + T[3,0,1])/2:
}
\vskip 4mm

We remind the reader that $S_1, S_2$ are \emph{relative} invariants
under contact transformations; they transform by \eqref{S1S2trans} and 
\eqref{S1S2trans2}.\\

\noindent{\bf II. Euler-Lagrange Examples Revisited.}

Given a hyperbolic Monge-Amp\`ere PDE \eqref{MAequation}, 
the code above can be used to decide whether it is Euler-Lagrange and, 
when it \emph{is}, its Euler-Lagrange type (positive, negative, degenerate). 

For the examples \ref{K=-1ex} and \ref{K=1ex} in Section \ref{InvMA}, 
we remind the reader of their classical PDE forms; then,
we continue the calculation above by computing $S_2$ and $\det(S_1)$
for these examples as well as for the Monge-Amp\`ere equation satisfying
\eqref{ABCDE}. 

\begin{enumerate}[{\bf 1.}]
\item{{\bf $K = -1$ Surfaces in $\E^3$.} (Example \ref{K=-1ex}, Section \ref{InvMA})

Suppose that a surface $S\subset \E^3$ is a graph with the position vector $\x = (x,y,z(x,y))$.
Letting $p:=z_x$ and $q:=z_y$, we have that, pulled back to $S$,
\[
		\theta:=\ed z - p\ed x - q\ed y = 0.
\]
The oriented unit normal along $S$ is just
\[
	\e_3 = (-p, -q, 1)^T(1+p^2+q^2)^{-1/2}.
\]
Direct calculation yields the area form
\[
	\ed\mathcal{A} =\frac{1}{2}(\ed \x\times\ed \x)\cdot \e_3 = (1+p^2+q^2)^{1/2}\ed x\W\ed y
\]
and
\[
	K\ed\mathcal{A} =\frac{1}{2}(\ed \e_3\times\ed \e_3)\cdot \e_3 = (1+p^2+q^2)^{-3/2} \ed p\W\ed q,
\]
where $K$ is the Gauss curvature of $S$.

It follows that the system for $K = -1$ is characterized by
\[
	\theta =0, \qquad \ed p\W \ed q + (1+p^2+q^2)^2 \ed x \W\ed y = 0.
\]
Noting that $\ed x\W\ed y \ne 0$ on $S$, the PDE form of this system is:
\begin{equation}\label{K=-1PDE}
	(z_{xx}z_{yy}- z_{xy}^2) + (1+z_x^2+z_y^2)^2 = 0. 
\end{equation}}

\item{{\bf $K = 1$ Time-like Surfaces in $\E^{2,1}$.} (Example \ref{K=1ex}, Section \ref{InvMA})

Since we use the signature $(+,-,+)$ for the Lorentzian metric on $\E^{2,1}$, the dot and cross 
products are defined, respectively, by
\[
	\u\cdot {\bf v} = u_1v_1 - u_2 v_2 + u_3 v_3,\qquad \u\times{\bf v} = \(\begin{array}{c}
								u_2v_3 - u_3v_2\\
								u_1v_3 - u_3v_1\\
								u_1v_2 - u_2v_1
						\end{array}\).
\]
Let $S$ be a time-like surface with the position vector $\x$ and 
a pseudo-orthonormal frame field
$(\e_1,\e_2,\e_3)$ attached to it ($\e_1,\e_2$ being tangent to $S$). 
The area form $\ed \mathcal{A}$ and the Gauss curvature $K$ can be computed by
	\[
		\ed \mathcal{A} = \frac{1}{2}(\ed \x \times \ed \x)\cdot \e_3,
	\]
and
	\[
		K\ed \mathcal{A} = \frac{1}{2}(\ed \e_3\times \ed \e_3)\cdot \e_3,
	\]
where the $2$-forms on the right-hand-sides are their pull-backs to $S$.
	
When $S$ is a graph $\x = (x,y,z(x,y))$, it is easy to find that 
	\[
		\e_3 = (-p, q, 1)^T(1+p^2 - q^2)^{-1/2}.
	\]
Direct calculation yields that
	\[
		\ed \mathcal{A} =(1+p^2 - q^2)^{1/2}\ed x \W\ed y , 
		\quad K\ed \mathcal{A} = -(1+p^2-q^2)^{-3/2}\ed p \W\ed q.
	\]
Therefore, the equation for $K = 1$ is
	\begin{equation}\label{K=1PDE}
		(z_{xx}z_{yy} - z_{xy}^2) + (1+ z_x^2 - z_y^2)^2 = 0.
	\end{equation}
}
\end{enumerate}

Now, calculating using Maple\texttrademark, we find:\footnote{Note that, in
\eqref{ABCDE}, $2p^2z^2+2p^2+zq>0$, by hyperbolicity.}
{\small
\[
\begin{array}{|c||c|c|c|}
	\hline
	\text{Eq.} &  S_2 & \det(S_1)&\text{E-L~Type}\\
	\hline\hline
	\eqref{K=-1PDE}&\bs 0&-\dfrac{1}{16}(p^2+q^2+1)& -\\[0.5em]
	\eqref{K=1PDE}&\bs 0&\dfrac{1}{16}(p^2 - q^2+1)& +\\[0.5em]
	\eqref{ABCDE}&\bs 0&\dfrac{z^2q^4(4p^2z^5-16p^2z^3+qz^4-20p^2z-8qz^2-q)^2}{32(2p^2z^2+2p^2+zq)^3(z^2+1)^6}& +\\[0.5em]
	\hline
\end{array}
\]}

The Euler-Lagrange types are consistent with our earlier observation when calculation was made using
differential forms. Here, we have used the following code for 
computing $\det(S_1)$ and $S_2$ associated to \eqref{ABCDE}, which can be 
easily modified to work
for \eqref{K=-1PDE} and \eqref{K=1PDE}.
\\

{\fontfamily{pcr}\selectfont \scriptsize
\setlength\parindent{0pt}

> ABCDE:= \{eA = 2*q*z*(z$\W$2+1)$\W$3, 

\hspace{5mm}eB = 2*q$\W$2*(z$\W$2+1)$\W$2*(4*p$\W$2*z$\W$3 - q*z$\W$2 + 4*p$\W$2*z + 3*q),

\hspace{5mm}eC = -2*p*q*(z$\W$2+1)$\W$3*(4*p$\W$2*z+q),

\hspace{5mm}eD = (z$\W$2+1)*(4*p$\W$2*z$\W$3+q*z$\W$2+4*p$\W$2*z - q)*(2*p$\W$2*z$\W$2+q*z+2*p$\W$2),

\hspace{5mm}eE = -q$\W$3*(4*p$\W$2*z$\W$5+q*z$\W$4 - 16*p$\W$2*z$\W$3 - 8*q*z$\W$2 - 20*p$\W$2*z - q)\}:

> BCDM:=  subs(ABCDE, \{B = eB/eE, C = eC/eE, D = eD/eE, 

\hspace{\codesp} M = sqrt((eA*eE + eC$\W$2 - eB*eD))/eE\}):

> dBCDM:= \{\}: ddBCDM:= \{\}:

> for t in [B,C,D,M] do:

\hspace{\codesp}  for i from 1 to 5 do:
  
\hspace{\codesp}      dBCDM:= \{op(dBCDM), t[i] = diff(subs(BCDM, t), var[i])\}:
      
\hspace{\codesp} od:

\hspace{\codesp}od:

> for t in [B,C,D,M] do:

\hspace{\codesp}    for i from 1 to 5 do:
    
\hspace{\codesp}        for j from 1 to 5 do:
        
 \hspace{\codesp}          ddBCDM:= \{op(ddBCDM), t[i,j] = diff(diff(subs(BCDM, t), var[i]), var[j])\}:
 
 \hspace{\codesp}   od:

\hspace{\codesp} od:

\hspace{\codesp}od:

> for i from 5 to 8 do:

\hspace{\codesp}   simplify(subs(BCDM union dBCDM union ddBCDM, V[i]));

\hspace{\codesp}od;

> subs(BCDM union dBCDM union ddBCDM, V[1]*V[4]-V[2]*V[3]):simplify(\%)
\\

}

\noindent{\bf III. The Lagrangian.}

Classically, a central object in a 2D variational problem is
the functional:
\[
		\int_{\mathcal{D}} L(x,y,z,z_x,z_y)\ed x\ed y,
\]
where $L$ is the Lagrangian (function) and $\mathcal{D}\subset \R^2$ a fixed compact domain.

The Euler-Lagrange equation associated to a fixed-boundary variation is
\begin{equation}\label{classicalEL}
	L_z  =  \frac{\ed L_{z_x}}{\ed x}  + \frac{\ed L_{z_y}}{\ed y}.
\end{equation}

A coordinate-independent formulation (see \cite[Section 1.2]{BGG}) considers instead
a fixed boundary variation in $J^1(\R^2,\R)$ (with standard coordinates
$(x,y,z,p,q)$) by Legendre surfaces\footnote{That is, surfaces that annihilate the 
pull-back of $\ed z - p\ed x -q\ed y$.} $S$ and the functional:
\[
	\int_S \Lambda,
\]
where $\Lambda = L(x,y,z,p,q)\ed x\W\ed y$.
This variational problem may seem less restrictive than the classical one,\footnote{Note that a fixed-boundary
variation in the classical formulation may not lift to be a fixed-boundary variation
by Legendre surfaces in $J^1(\R^2,\R)$.} but the corresponding Euler-Lagrange equation is the same as \eqref{classicalEL},
which is Monge-Amp\`ere 
 \eqref{MAequation} with
\[
	A = 0, \quad B = L_{pp}, \quad C = L_{pq}, \quad D = L_{qq}, 
\]
\[
	E = L_{xp} + p L_{zp} + L_{yq} + qL_{zq} - L_z.
\]

More generally but still in the 2D case, $\Lambda$ can be any
$2$-form on $J^1(\R^2,\R)$.
The `stationary point' of a corresponding fixed-boundary variational problem
is characterized by a more general Euler-Lagrange system (see \cite[Section 1.2.2]{BGG}). It is this latter notion 
that we have in mind when we refer to Euler-Lagrange systems.

Theoretically, the inverse problem: 
\begin{quote}
\emph{Given a hyperbolic Monge-Amp\`ere
system with $S_2 = \bs 0$, can one construct an associated Lagrangian $2$-form $\Lambda$ in some local coordinates?}
\end{quote}
is solved in \cite[Theorems 1.2 and 2.2]{BGG}.
Here we only outline the steps
of a construction and present some examples.

Start with a hyperbolic Monge-Amp\`ere system with $S_2 = \bs0$.
\begin{enumerate}[{\bf Step 1.}]
	\item{Find a $1$-adapted coframing $\bs\omega = (\omega^0,\ldots,\omega^4)$;}
	\item{Compute the $1$-form $\phi_0$ in \eqref{StrEqnMA}, which is determined, using
		$\bs\omega$;}
	\item{Since, $S_2 = \bs 0$, $\phi_0$ must be closed; thus, find a function 
			$\lambda$ such that 
			$\ed \lambda = 2\lambda \phi_0$;}
	\item{It would then follow that $\Pi := \lambda\omega^0\W(\omega^1\W\omega^2 - \omega^3\W\omega^4)$ is closed; hence, it is a
	Poincar\'e-Cartan form.
			Now find\footnote{This step can involve using the proof of the Poincar\'e Lemma.} a $2$-form $\Lambda$ such that 
			$\ed \Lambda = \Pi$. This $\Lambda$ is a desired Lagrangian.
}			
\end{enumerate}

Now, for \eqref{K=-1PDE}, \eqref{K=1PDE} and the Monge-Amp\`ere equation satisfying
\eqref{ABCDE}, we follow these steps to compute their representative Lagrangian 2-forms and
summarize intermediate expressions in the following tables. 

{\footnotesize
\[
\begin{array}{|c||c|}
	\hline
	\text{Eq.} &  \eqref{K=-1PDE}\\
	\hline\hline
	\phi_0 &\dfrac{p\ed p + q\ed q}{1+p^2+q^2}\\
	\lambda& 1+p^2+q^2\\
	\Pi& 4\ed x\W\ed y\W\ed z + \dfrac{4(\ed z-p\ed x-q\ed y)}{(1+p^2+q^2)^2}\ed p\W\ed q\\
	\Lambda&4z~ \ed x \W \ed y + \dfrac{4z}{(1+p^2+q^2)^2}\ed p\W\ed q
				 - \dfrac{2}{1+p^2+q^2}(\ed x\W\ed q- \ed y\W\ed p)\\[0.5em]
	\hline
\end{array}
\]	

\[
\begin{array}{|c||c|}
	\hline
	\text{Eq.} &  \eqref{K=1PDE}\\
	\hline\hline
	\phi_0 & \dfrac{p\ed p - q\ed q}{1+p^2 -q^2}\\
	\lambda& 1+ p^2 - q^2 \\
	\Pi& 4\ed x\W\ed y\W\ed z +\dfrac{4(\ed z - p\ed x - q\ed y) }{(1+p^2-q^2)^2}\ed p\W\ed q\\
	\Lambda&4z~\ed x \W\ed y + \dfrac{4z}{(1+p^2 - q^2)^2}\ed p\W\ed q
		-\dfrac{2}{1+p^2 - q^2}(\ed x \W\ed q + \ed y \W\ed p)\\[0.5em]
	\hline
\end{array}
\]	}

For the Monge-Amp\`ere equation satisfying \eqref{ABCDE}, we record only $\lambda$
and $\Lambda$, since the expressions of $\phi_0$ and $\Pi$ are rather complicated.

{\footnotesize
\[\def\arraystretch{2}
\begin{array}{|c||c|}
	\hline
	\text{Eq.} &  \text{Monge-Amp\`ere Equation }\eqref{MAequation} \text{ satisfying } \eqref{ABCDE}\\
	\hline\hline
	\lambda& \dfrac{q^2(4p^2z^5 - 16p^2z^3 + qz^4 - 20 p^2z - 8qz^2-q)^2}{(2p^2z^2+2p^2+qz)^2(z^2+1)^4}\\[0.5em]
	\hline\hline{\rule{0pt}{6.8ex}}
	\dfrac{(2z^2p^2+2p^2+qz)^{1/2}}{8\sqrt{2}}\Lambda&
	\begin{split}	
				&	4(z^2+1)\ed y\W\ed p + \dfrac{(px+2qy)(8p^2z+2q)(z^2+1)}{q(2z^2p^2+2p^2+qz)}\ed z\W\ed p\\
				& - \dfrac{z}{q}\ed x\W\ed q  + \dfrac{2q(z^2-1)}{z^2+1}\ed x\W\ed y
				+\Sigma~ \ed z\W\ed q
	\end{split}				
					\\[2em]
	\hline
\end{array}
\]	
}
In the above,
\begin{align*}
	\Sigma &= \frac{1}{2(z^2+1)q^2(2p^2z^2+2p^2+zq)}\Big(-8p^3z^6 + (-16p^4x - 16p^3qy - 4pq)z^5\\
	&+(-12p^2qx - 4pq^2y - 24p^3)z^4+(-32p^4x-32p^3qy-3q^2x - 8pq)z^3\\
	&+(-12p^2qx - 8pq^2y - 24p^3)z^2 + (-16p^4x - 16p^3qy+q^2x - 4pq)z\\
	& -4pq^2y - 8p^3\Big).
\end{align*}

Finally, we remind the reader that $\Lambda$ is not unique in the sense that one can
add to it any exact $2$-form and any $2$-form in the contact ideal generated by
$\ed z - p\ed x - q\ed y$. For details, see \cite[Section 1.1]{BGG}.

\end{appendix}

\bibliographystyle{alpha}
\normalbaselines 

\newcommand{\etalchar}[1]{$^{#1}$}

\end{document}